\documentclass[10pt]{article}

\usepackage[margin=1in]{geometry}

\usepackage{amsmath}
\usepackage{amsthm}
\usepackage{graphicx}
\usepackage{color}
\usepackage{amsfonts}
\usepackage{amssymb}
\usepackage{psfrag}
\usepackage{wrapfig}
\usepackage{amscd}
\usepackage{url}

\newcommand{\re}{\mathbb{R}}
\newcommand{\cpx}{\mathbb{C}}

\newcommand{\N}{\mathbb{N}}

\newcommand{\half}{\frac{1}{2}}
\newcommand{\lmd}{\lambda}

\newcommand{\nn}{\nonumber}

\def\af{\alpha}
\def\bt{\beta}

\newcommand{\sig}{\sigma}

\newcommand{\reff}[1]{(\ref{#1})}
\newcommand{\pt}{\partial}
\newcommand{\prm}{\prime}

\newcommand{\mc}[1]{\mathcal{#1}}

\newcommand{\bdes}{\begin{description}}
\newcommand{\edes}{\end{description}}

\newcommand{\bal}{\begin{align}}
\newcommand{\eal}{\end{align}}

\newcommand{\bnum}{\begin{enumerate}}
\newcommand{\enum}{\end{enumerate}}

\newcommand{\bit}{\begin{itemize}}
\newcommand{\eit}{\end{itemize}}

\newcommand{\bea}{\begin{eqnarray}}
\newcommand{\eea}{\end{eqnarray}}
\newcommand{\be}{\begin{equation}}
\newcommand{\ee}{\end{equation}}

\newcommand{\baray}{\begin{array}}
\newcommand{\earay}{\end{array}}

\newcommand{\bsry}{\begin{subarray}}
\newcommand{\esry}{\end{subarray}}

\newcommand{\bca}{\begin{cases}}
\newcommand{\eca}{\end{cases}}

\newcommand{\bcen}{\begin{center}}
\newcommand{\ecen}{\end{center}}

\newcommand{\bbm}{\begin{bmatrix}}
\newcommand{\ebm}{\end{bmatrix}}

\newcommand{\bmx}{\begin{matrix}}
\newcommand{\emx}{\end{matrix}}

\newcommand{\bpm}{\begin{pmatrix}}
\newcommand{\epm}{\end{pmatrix}}

\newcommand{\btab}{\begin{tabular}}
\newcommand{\etab}{\end{tabular}}

\newcommand{\thmlist}{
\begin{list}{Step 1}
{\setlength{\leftmargin}{0.6 in}\setlength{\labelwidth}{0.5 in}} }
\newcommand{\alglist}{
\begin{list}{Step 1}
{\setlength{\leftmargin}{1.1 in}\setlength{\labelwidth}{1.0 in}} }
\theoremstyle{plain}
\newtheorem{theorem}{Theorem}[section]
\newtheorem{pro}[theorem]{Proposition}

\newtheorem{lemma}[theorem]{Lemma}

\theoremstyle{definition}
\newtheorem{example}[theorem]{Example}

\renewcommand{\subsection}[1]{
    \stepcounter{subsection}
    \settowidth{\hangindent}{\bf\thesubsection.~}
    \hangafter=1
    \bigskip\noindent
    {\bf\hbox{\thesubsection.~}#1}\par
    \nobreak
    \medskip
}
\begin{document}

\title{ First Order Conditions for Semidefinite Representations of Convex Sets Defined by
Rational or Singular Polynomials
\author{Jiawang Nie\footnote{Department of Mathematics,
University of California at San Diego, 9500 Gilman Drive, 
La Jolla, CA 92093, USA. Email: njw@math.ucsd.edu.}
}
\date{June 18, 2008}
}

\maketitle

\begin{abstract}
A set is called {\it semidefinite representable} or {\it semidefinite programming (SDP) representable} 
if it can be represented
as the projection of a higher dimensional set which is represented 
by some {\it Linear Matrix Inequality (LMI)}.
This paper discuss the semidefinite representability conditions
for convex sets of the form 
$S_{\mc{D}}(f) =\{ x\in \mc{D}:  f(x) \geq 0 \}$.
Here $\mc{D}=\{x\in\re^n:g_1(x)\geq 0, \cdots, g_m(x) \geq 0 \}$ is a convex domain
defined by some ``nice" concave polynomials $g_i(x)$
(they satisfy certain concavity certificates), 
and $f(x)$ is a polynomial or rational function.
When $f(x)$ is concave over $\mc{D}$, 
we prove that $S_{\mc{D}}(f) $ has some explicit semidefinite representations 
under certain conditions called {\it preordering concavity} or {\it q-module concavity},
which are based on the Positivstellensatz certificates
for the first order concavity criteria:
\[
f(u) + \nabla f(u)^T(x-u) -f(x) \geq 0, \quad \forall \, x, u \in \mc{D}.
\]
When $f(x)$ is a polynomial or rational function having singularities on the boundary of $S_{\mc{D}}(f)$,
a {\it perspective transformation} is introduced 
to find some explicit semidefinite representations for $S_{\mc{D}}(f)$
under certain conditions.
In the particular case $n=2$, if the Laurent expansion of $f(x)$ around one singular point
has only two consecutive homogeneous parts, 
we show that $S_{\mc{D}}(f)$ always admits 
an explicitly constructible semidefinite representation.
\end{abstract}
 
 
\section{Introduction}

Semidefinite programming (SDP) \cite{BTN, NN94, N06, SDPhandbook}
is an important convex optimization problem.
It has wide applications in combinatorial optimization, control theory
and nonconvex polynomial optimization as well as many other areas.
There are efficient numerical algorithms and standard packages
for solving semidefinite programming.
Hence, a fundamental problem in optimization theory is
what sets can be presented by semidefinite programming.
This paper discusses this problem.

A set $S$ is said to be {\it Linear Matrix Inequality (LMI) representable} if
\be \nn
S=\{x\in \re^n: A_0+A_1 x_1 + \cdots + A_nx_n \succeq 0\}
\ee
for some symmetric matrices $A_i$. Here the notation $X
\succeq 0 \,(\succ 0)$ means $X$ is positive
semidefinite (definite).  
The above is then called an LMI representation for $S$.
If $S$ is representable as the projection of
\be \nn
\hat S = \left\{ (x, u) \in \re^{(n+N)}: A_0 +
\sum_{i=1}^n A_i x_i + \sum_{j=1}^N B_j u_j \succeq 0 \right\}
\subset \re^{(n+N)},
\ee
that is,
$S = \left\{ x \in \re^{n}: \exists u \in \re^n, \
(x,u) \in \hat S \right\}
$,
for some symmetric matrices $A_i$ and $B_j$,
then $S$ is called {\it semidefinite representable} or {\it semidefinite programming (SDP) representable}.
The lifted LMI above is then called a
{\it semidefinite representation}, {\it SDP representation}
or {\it lifted LMI representation} for $S$.
Sometimes, we also say $S=\hat S$ if $S$ equals the projection of the lift $\hat S$.

Nesterov and Nemirovski
(\cite{NN94}), Ben-Tal and Nemirovski (\cite{BTN}), and Nemirovsky (\cite{N06})
gave collections of examples of SDP representable sets.
Thereby leading to the fundamental question which sets are SDP representable?
%In his ICM 2006 plenary lecture \cite{N06}, Nemirovsky remarked that
%this question seems to be completely open. 
Obviously, to be SDP representable, $S$ must be convex and semialgebraic.
Is this necessary condition also sufficient?
What are the sufficient conditions for $S$ to be SDP representable?
Note that not every convex semialgebraic set is LMI representable
(see Helton and Vinnikov \cite{HV}).

\medskip
\noindent
{\bf \large Prior work}
When $S$ is a convex set of the form
$\{x\in \re^n:\, g_1(x)\geq 0, \cdots, g_m(x)\geq 0\}$
defined by polynomials $g_i(x)$,
there is recent work on the SDP representability of $S$.
%Parrilo \cite{Par06} and Lasserre \cite{Las06} independently
%proposed a natural
%construction of lifted LMIs using moments and sum of squares techniques
%with the aim  of producing  SDP representations.
Parrilo \cite{Par06} gave a construction of lifted LMIs 
using moments and sum of squares techniques,
and proved the construction gives an SDP representation 
in the two dimensional case when the boundary
of $S$ is a single rational planar curve of genus zero.
Lasserre \cite{Las06} showed the construction can give arbitrarily accurate
approximations to compact $S$,
and the construction gives a lifted LMI for $S$
under some algebraic properties called {\it S-BDR} or {\it PP-BDR}, i.e.,
requiring {\it almost all} positive affine polynomials on $S$
have certain SOS representations with uniformly bounded degrees.
Helton and Nie \cite{HN1} proved that 
the convex sets of the form
$\{x\in \re^n:\, g_1(x)\geq 0, \cdots, g_m(x)\geq 0\}$
are SDP representable if every $g_i(x)$ is sos-concave
($-\nabla^2 g_i(x) = G_i(x)^T G_i(x)$ for some possibly nonsquare matrix polynomial $G_i(x)$),
or every $g_i(x)$ is strictly quasi-concave on $S$,
or a mixture of the both.
Later, based on the work \cite{HN1},
Helton and Nie \cite{HN2} proved a very general result that 
a compact convex semialgebraic set $S$ is always SDP representable 
if the boundary of $S$ is nonsingular and has positive curvature.
This sufficient condition is not far away from being necessary:
the boundary of a convex set has nonnegative curvature when it is nonsingular.
So the only unaddressed cases for SDP representability are that
the boundary of a convex set has zero curvature somewhere or
has some singularities.

\medskip \noindent
{\bf \large Contributions} \,
The results in \cite{HN1, HN2, Las06} are more on the theoretical existence
of SDP representations. The constructions given there 
might be too complicated to be useful for computational purposes.
And these results sometimes need check conditions of 
Hessians of defining polynomials, 
which sometimes are difficult or inconvenient to verify in practice.
However, in many applications, people often want explicit and simple 
semidefinite representations.   
Thus some ``simple" SDP representations and conditions justifying them
are favorable in practical applications. 
All these practical issues motivate this paper. 
Our contributions come in the following three aspects.

%This paper discuss the SDP representation of convex sets in $\re^n$
%defined by polynomial or rational functions.
%Let
%\[
%S(f_1,\cdots,f_k) = \{x \in\re^n: f_1(x) \geq 0,\cdots,f_k(x)\geq 0\}
%\]
%be a semialgebraic set defined by functions $f_1,\cdots,f_k$.
%When every $f_i$ is an sos-concave polynomial, it was shown that
%$S(f_1,\cdots,f_k)$ is SDP representable
%and the representation can be constructed explicitly.
%When $S(f_1,\cdots,f_k)$ is compact and its boundary has positive curvature,
%it was shown it also allows an SDP representation.
%However, the semidefinite representability
%for some particular convex sets
%is not addressed in the prior work

\medskip

{\it First}, there are some convex sets defined by polynomials
that are not concave in the whole space $\re^n$
but concave over a convex domain $\mc{D}\subset \re^n$.
For instance, for convex set
$\{x\in\re^2: x_2-x_1^3\geq 0, x_1\geq 0\}$,
the defining polynomial $x_2-x_1^3$ is not  concave when $x_1 <0$,
but is concave over the domain $\re_+\times \re$.
However, this set allows an SDP representation, e.g., 
\[
\left\{(x_1,x_2): \, \exists \, u, \,\,
\bbm x_1 & u \\ u & x_2 \ebm \succeq 0, \,
\bbm 1 &x_1  \\ x_1 & u \ebm \succeq 0
\right\}.
\]
For convex sets given in the form
$S_{\mc{D}}(f)=\{x \in \mc{D}: f(x)\geq 0\}$,
where $f(x)$ is a polynomial concave over a convex domain $\mc{D}$,
we prove some sufficient conditions for 
semidefinite representability of $S_{\mc{D}}(f)$
and give explicit SDP representations.
This will be discussed in Section~\ref{sec:poly}.

\medskip
{\it Second}, there are some convex sets defined
by rational functions (also called rational polynomials) 
which are concave over a convex domain $\mc{D}$ of $\re^n$.
If we redefine them by using polynomials, the concavity of rational functions
might not be preserved.
For instance, the unbounded convex set
\[
\left\{x\in \re_+^2:\, 1 - \frac{1}{x_1x_2} \geq 0 \right\}
\]
is defined by a rational function concave over $\re_+^2$
($\re_+$ is the set of nonnegative real numbers).
This set can be equivalently defined by polynomials
\[
\left\{x\in \re^2:\, x_1x_2 - 1 \geq 0, x_1\geq 0, x_2 \geq 0 \right\}.
\]
But $x_1x_2-1$ is not concave anywhere.
The prior results in \cite{HN1,HN2} do not imply
the SDP representability of this set.
However, this set is SDP representable, e.g.,
\[
\left\{x\in \re^2: \bbm x_1 & 1 \\ 1 & x_2 \ebm \succeq  0 \right\}.
\]
For convex sets given in the form
$S_{\mc{D}}(f)=\{x \in \mc{D}: f(x)\geq 0\}$,
where $f(x)$ is a rational function concave over a convex domain $\mc{D}$,
we prove some sufficient conditions for 
semidefinite representability of $S_{\mc{D}}(f)$
and give explicit SDP representations.
This will be discussed in Section~\ref{sec:RatCvx}.

\medskip
{\it Third}, there are some convex sets that are defined
by polynomials or rational functions which are singular on the boundary.
For instance, the set
\[
\{x\in \re^2:\, x_1^2-x_1^3-x_2^2 \geq 0, x_1 \geq 0 \}
\]
is convex, and the origin is on the boundary.
The polynomial $x_1^2-x_1^3-x_2^2$ is singular at the origin,
i.e., its gradient vanishes at the origin.
The earlier results in \cite{HN1,HN2} do not imply the SDP representability of this set.
However, this set can be equivalently defined as
\[
\left\{(x_1,x_2) \in \re_{+} \times \re:\,
x_1- x_1^2-\frac{x_2^2}{x_1} \geq 0 \right\},
\]
a convex set defined by a concave rational function over the domain
$\re_{+} \times \re$. 
By Schur's complement, we know it can be represented as
\[
\left\{(x_1,x_2):\,
\bbm x_1 & x_2 & x_1 \\ x_2 & x_1 & 0 \\ x_1 & 0 & 1  \ebm \succeq 0
\right\}.
\]
It is an LMI representation without projections.
The technique of Schur's complement works only for 
very special concave rational functions,
and is usually difficult to be applied for general cases.
%On the other hand, note that
%a convex set with singularities can be defined by different
%rational functions
%and some of them might not be concave.
%For instance, this convex set can also be defined as
%\[
%\left\{(x_1,x_2) \in \re_{+} \times \re:\,
%1- x_1-\frac{x_2^2}{x_1^2} \geq 0 \right\}.
%\]
%But its defining rational function is not concave over the domain $\re_+ \times \re$.
For singular convex sets of the form
$S_{\mc{D}}(f)=\{x \in \mc{D}: f(x)\geq 0\}$,
where $f(x)$ is a polynomial or rational function with singularities on the boundary,
we give some sufficient conditions for 
semidefinite representability of $S_{\mc{D}}(f)$
and give explicit SDP representations.
In the particular case $n=2$, we show that 
$S_{\mc{D}}(f)$ always admits an explicitly constructible SDP representation
when the Laurent expansion of $f(x)$ around one singular point 
has only two consecutive homogeneous parts.
This will be discussed in Section~\ref{sec:sig}.

\medskip

In this paper, we always assume
$\mc{D}=\{x\in\re^n:g_1(x)\geq 0, \cdots, g_m(x) \geq 0 \}$ 
is a convex domain defined by some nice concave polynomials $g_i(x)$.
Here ``nice" means that they satisfy certain concavity certificates.
For instance, a very useful case is $\mc{D}$ is a polyhedra.
We do not require $\mc{D}$ or $S_{\mc{D}}(f)$ to be compact,
as required by \cite{HN1,HN2,Las06}.
When $f(x)$ is concave over $\mc{D}$,
the sufficient conditions for SDP representability of $S_{\mc{D}}(f)$
proven in this paper
are based on some certificates for the first order concavity criteria:
\[
f(u) + \nabla f(u)^T(x-u) -f(x) \geq 0, \quad \forall \, x, u \in \mc{D}.
\]
Some Positivstellensatz certificates
like Putinar's Positivstellensatz \cite{Putinar} or Schm\"udgen's Positivstellensatz \cite{Smg}
for the above
can be applied to justify some explicitly constructible 
SDP representations for $S_{\mc{D}}(f)$.

Throughout this paper, $\re$ (resp. $\N$) denotes the set of real numbers (resp. nonnegative integers).
For $\af \in \N^n$ and $x\in \re^n$, denote $|\af| = \af_1+\cdots+\af_n$
and $x^\af = x_1^{\af_1}\cdots x_n^{\af_n}$.
$B(u,r)$ denotes the ball $\{x\in\re^n: \|x-u\|_2 \leq r\}$.
A vector $x\geq 0$ means all its entries are nonnegative.
A polynomial $p(x)$ is said to be a sum of squares or sos 
if there finitely many polynomials $q_i(x)$ such that 
$p(x) = \sum q_i(x)^2$. 
A matrix polynomial $H(x)$ is called a sum of squares or sos 
if there is a matrix polynomial $G(x)$ such that $H(x) = G(x)^TG(x)$.

\section{Convex sets defined by polynomials concave over domains}
\label{sec:poly}

In this section, consider the convex set
$S_{\mc{D}}(f) = \{x \in \mc{D}:\, f(x)\geq 0\}$
defined by a polynomial $f(x)$.
Here $\mc{D}=\{x \in \re^n: g_1(x)\geq 0, \cdots, g_m(x) \geq 0 \}$ is a convex domain.
When $f(x)$ is concave on $\mc{D}$,
it must hold
\[
-R_f(x,u):=f(u) + \nabla f(u)^T(x-u) -f(x) \geq 0, \, \forall\, x, u \in \mc{D}.
\]
The difference $R_f(x,u)$ is the first order Lagrange remainder.

\subsection{q-module convexity and preordering convexity}

Now we introduce some types of definitions about convexity/concavity.
Define $g_0(x)=1$.
We say $f(x)$ is {\it q-module convex over $\mc{D}$}
if it holds
\[
R_f(x,u) = \sum_{i=0}^m g_i(x) \left( \sum_{j=0}^m g_j(u) \sig_{ij}(x,u)  \right)
\]
for some sos polynomials $\sig_{ij}(x,u)$.
Then define $f(x)$ to be {\it q-module concave over $\mc{D}$}
if $-f(x)$ is  q-module convex over $\mc{D}$.
We say $f(x)$ is {\it preordering convex over $\mc{D}$}
if it holds
\[
R_f(x,u) = \sum_{\nu\in \{0,1\}^m } g_1^{\nu_1}(x)\cdots g_m^{\nu_m}(x)
\left(  \sum_{\mu\in \{0,1\}^m } g_1^{\mu_1}(u)\cdots g_m^{\mu_m}(u)  \sig_{\nu,\mu}(x,u)  \right)
\]
for some sos polynomials $\sig_{\nu,\mu}(x,u)$.
Similarly, $f(x)$ is called {\it preordering concave over $\mc{D}$}
if $-f(x)$ is preordering convex over $\mc{D}$.
Obviously, the q-module convexity implies preordering convexity,
which then implies the convexity,
but the converse might not be true.

We remark that the defining polynomials $g_i(x)$ are not unique for the domain $\mc{D}$.
When we say $f(x)$ is  q-module or preordering convex/concave over $\mc{D}$,
we usually assume a certain set of defining polynomials $g_i(x)$ is clear in the context.

In the special case $\mc{D}=\re^n$,
the definitions of q-module convexity and preordering convexity
coincide each other, and then are specially called
{\it first order sos convexity}.
And {\it first order sos concavity} is defined in a similar way.
Recall that a polynomial $f(x)$ is {\it sos-convex} if
its Hessian $\nabla^2 f(x)$ is sos (see \cite{HN1}).
An interesting fact is if $f(x)$ is sos-convex 
then it must also be first order sos convex.
This is due to that
\begin{align*}
& \quad f(x) - f(u) - \nabla(u)^T(x-u)  \\
&= (x-u)^T \left(\int_0^1\int_0^t f^{\prm\prm}(u+s(x-u)) \, ds \, dt \right) (x-u) \\
& = (x-u)^T \left(\int_0^1\int_0^t F(u+s(x-u))^T F(u+s(x-u)) \, ds \, dt \right) (x-u)
\end{align*}
is an sos polynomial (see Lemma~3.1 of \cite{HN1}).

%
%
%\begin{pro}
%Suppose $f(x)$ is convex on $\mc{D}$.
%If the Hessian $\nabla^2 f(x)$ belongs to
%the quadratic module (resp. preporder)
%generated by $g_1(x), \cdots, g_m(x)$,
%then $(x)$ is q-module convex (resp. preporder convex)
%over $\mc{D}$.
%\end{pro}
%
%
%
%\begin{pro}
%Suppose $f(x)$ has a minimizer $\xi$ in the interior of $\mc{D}$.
%Then we have
%\bit
%\item [(i)] If $f(x)$ is q-module convex over $\mc{D}$,
%then $f(x)-f(\xi)$ belongs to the quadratic module generated by $g_1,\cdots,g_m$.
%
%\item [(ii)] If $f(x)$ is preordering convex over $\mc{D}$,
%then $f(x)-f(\xi)$ belongs to the preordering generated by $g_1,\cdots,g_m$.
%
%\item [(iii)] When $\mc{D}=\re^n$, if a polynomial $f(x)$ is first order sos convex,
%then $f(x)-f(\xi)$ must be sos.
%
%\eit
%\end{pro}
%

\begin{example} \label{exm:2dot1}
The bivariate polynomial $f(x)=x_1^3+x_1^2x_2+x_1x_2^2+x_2^3$
is convex over the nonnegative orthant $\re_+^2$.
It is also q-module convex with respect to $\re_+^2$.
This is due to the identity
\begin{align*}
R_f(x,u) =
\Big(\frac 13 u_1 + \frac 16 x_1\Big)  \Big( 4(x_1-u_1)^2 + 2 (x_1+x_2-u_1-u_2)^2\Big) + \\
\qquad \Big(\frac 13 u_2 + \frac 16 x_2\Big)  \Big( 4(x_2-u_2)^2 + 2 (x_1+x_2-u_1-u_2)^2\Big).
\end{align*}
\end{example}

\subsection{SDP representations}

Throughout this subsection, we assume the polynomials $f(x)$
and $g_i(x)$ are either all q-module concave or all preordering concave over $\mc{D}$.
For any $h(x)$ from the set $\{f(x), g_1(x), \cdots, g_m(x)\}$, we thus have either 
\[
-R_h(x,u) = \sum_{i=0}^m g_i(x) \left( \sum_{j=0}^m g_j(u) \sig_{ij}^{h}(x,u)  \right)
\]
for some sos polynomials $\sig_{ij}^{h}(x,u)$, or 
\[
-R_h(x,u) = \sum_{\nu\in \{0,1\}^m } g_1^{\nu_1}(x)\cdots g_m^{\nu_m}(x)
\left(  \sum_{\mu\in \{0,1\}^m } g_1^{\mu_1}(u)\cdots g_m^{\mu_m}(u)  \sig_{\nu,\mu}^{h}(x,u)  \right)
\]
for some sos polynomials $\sig_{\nu,\mu}^{h}(x,u)$. Let
$d_i = \lceil \half \deg(g_i) \rceil$,
$d_\nu = \lceil \half \deg_x(g_1^{\nu_1}\cdots g_m^{\nu_m} ) \rceil$ and 
\be \label{P:dgfrm}
\baray{c}
d_{qmod}^{(P)} = \underset{h \in \{f, g_1, \cdots, g_m\}}{\max}  \quad
\underset{0\leq i,j \leq m}{\max} \quad
\lceil \half \deg_x(g_i \sig_{ij}^{h}) \rceil \\
d_{pre}^{(P)} = \underset{h \in \{f, g_1, \cdots, g_m\}}{\max}  \quad
\underset{\nu\in \{0,1\}^m, \mu\in \{0,1\}^m }{\max} \quad
\lceil \half \deg_x(g_1^{\nu_1}\cdots g_m^{\nu_m} \sig_{\nu,\mu}^{h}) \rceil
\earay
\ee
where $\deg_x(\cdot)$ denotes the degree of a polynomial in $x$.
Then define $d=d_{qmod}^{(P)}$ (resp. $d=d_{pre}^{(P)}$)
when $f(x)$ and all $g_i(x)$
are q-module (resp. preordering) concave over $\mc{D}$.

Define matrices $G_\af^{(i)}$ and $G_\af^{(\nu)}$ such that
\be  \label{P:mprod}
\baray{rcl}
g_i(x)\mathfrak{m}_{d-d_i}(x)\mathfrak{m}_{d-d_i}(x)^T &= &
 \underset{\af\in \N^n: |\af| \leq 2d }{\sum}
G_\af^{(i)} x^{\af}  \\
g_1^{\nu_1}(x)\cdots g_m^{\nu_m}(x) \mathfrak{m}_{d-d_\nu}(x)\mathfrak{m}_{d-d_\nu}(x)^T 
&= &  \underset{\af\in \N^n: |\af| \leq 2d }{\sum}
G_\af^{(\nu)} x^{\af}.
\earay
\ee
Here $\mathfrak{m}_{k}(x)$ is the vector of all monomials with degrees $\leq k$.
Let $y$ be a vector multi-indexed by integer vectors in $\N^n$. Then define
\begin{align*}
N_i(y) &= \underset{\af\in \N^n: |\af| \leq 2d }{\sum}
G_\af^{(i)} y_\af, \quad i=0,1,\ldots, m,   \\
N_\nu(y) &= \underset{\af\in \N^n: |\af| \leq 2d }{\sum}
G_\af^{(\nu)} y_\af, \quad \nu \in \{0,1\}^m.
\end{align*}
Suppose the polynomial $f(x)$ is given in the form
$f(x) = \sum_{\af \in \N^n:\, |\af| \leq 2d } \,\, f_\af\, x^\af$.
Then define vector $f$ such that
$
f^T y =   \sum_{\af \in \N^n:\, |\af| \leq 2d } \,\, f_\af\, y_\af.
$
Define two sets
\begin{align}
\mc{L}_{qmod}^{\mc{D}}(f) & =\Big\{y \in \re_{\af\in \N^n, |\af|\leq 2d}: \, y_{0} = 1, \,
f^Ty \geq 0,\,  N_i(y) \succeq 0, \, \forall \, 0\leq i\leq m \Big\} \label{poLqmod}  \\
\mc{L}_{pre}^{\mc{D}}(f) & =\Big\{y \in \re_{\af\in \N^n, |\af|\leq 2d}: \, y_{0} = 1, \,
f^Ty \geq 0,\,  N_\nu(y) \succeq 0, \, \forall \, \nu \in \{0,1\}^m \Big\} \label{poLpre}
\end{align}
via their LMI representations.
Then $S_{\mc{D}}(f)$ is contained in the image of
both $\mc{L}_{qmod}^{\mc{D}}(f)$ and $\mc{L}_{pre}^{\mc{D}}(f)$
under the projection map:
$\rho(y) = (y_{1 0 \ldots 0}, \,\,\, y_{0 1 0 \ldots 0}, \,\,\, \ldots, \,\,\, y_{0 0 \ldots 0 1})$,
because for any $x\in S_{\mc{D}}(f)$ we can choose $y_\af = x^\af$
such that $y \in \mc{L}_{qmod}^{\mc{D}}(f)$ and $ y\in \mc{L}_{pre}^{\mc{D}}(f)$.
We say $S_{\mc{D}}(f)$ equals $\mc{L}_{qmod}^{\mc{D}}(f)$ (resp. $\mc{L}_{pre}^{\mc{D}}(f)$)
if $S_{\mc{D}}(f)$ equals the image of $\mc{L}_{qmod}^{\mc{D}}(f)$ (resp. $\mc{L}_{pre}^{\mc{D}}(f)$)
under this projection.
Similarly, we say $x\in \mc{L}_{qmod}^{\mc{D}}(f)$ (resp. $x\in \mc{L}_{pre}^{\mc{D}}(f)$)
if there exists $y\in \mc{L}_{qmod}^{\mc{D}}(f)$ (resp. $y\in \mc{L}_{pre}^{\mc{D}}(f)$)
such that $x = \rho(y)$.

\begin{lemma} \label{lem:qplinrep}
Assume $S_\mc{D}(f)$ has nonempty interior.
Let $\{x\in\re^n: a^Tx = b\} $ be a supporting hyperplane of $S_\mc{D}(f)$ such that
$
a^Tx \geq b,  \,\,\,\forall \,\, x \in S_\mc{D}(f)
$
and $a^Tu = b$ for some point $u\in S_\mc{D}(f) $.
\bdes
\item [(i)] If $f(x)$ and every $g_i(x)$ are q-module concave over $\mc{D}$,
then it holds
\[
a^Tx - b - \lmd f(x) =  \sum_{i=0}^m  g_i(x) \sig_i(x)
\]
for some scalar $\lmd \geq 0$ and sos polynomials $\sig_i(x)$ 
such that $\deg(g_i\sig_i) \leq 2d_{qmod}^{(P)}$.
\item [(ii)]
If $f(x)$ and every $g_i(x)$ are preordering concave over $\mc{D}$,
then it holds
\[
a^Tx - b - \lmd f(x) =  \sum_{\nu \in \{0,1\}^m }  g_1^{\nu_1}(x) \cdots g_m^{\nu_m}(x) \sig_\nu(x).
\]
for some scalar $\lmd \geq 0$ and sos polynomials $\sig_\nu(x)$
such that $\deg(g_1^{\nu_1}\cdots g_m^{\nu_m} \sig_\nu) \leq 2d_{pre}^{(P)}$.
\edes
\end{lemma}
\begin{proof}
Since $S_\mc{D}(f)$ has nonempty interior
and the polynomials $f(x), g_1(x), \cdots, g_m(x)$ are all cocnave, 
the first order optimality condition holds
at $u$ for convex optimization problem ($u$ is a minimizer)
\[
\min_{x}  \quad a^Tx  \quad
\mbox{ subject to }  \quad  f(x)\geq 0, \,  g_i(x)\geq 0,  \, i=1, \ldots, m.
\]
Hence there exist Lagrange multipliers $\lmd\geq 0, \lmd_1 \geq 0,\lmd_m \geq 0$ such that
\begin{align*}
a =  \lmd \nabla f(u) + \sum_{i=1}^m \lmd_i \nabla g_i(u), \quad
\lmd f(u) = \lmd_1 g_1(u) = \cdots =\lmd_m g_m(u) = 0.
\end{align*}
Thus the Lagrange function 
$a^Tx-b-\lmd f(x) -\sum_{i=1}^m \lmd_i g_i(x)$
has representation
\[
\lmd \Big ( f(u) + \nabla f(u)^T(x-u)  - f(x) \Big) +
\sum_{i=1}^m \lmd_i \Big ( g_i(u) + \nabla g_i(u)^T(x-u)  - g_i(x) \Big) .
\]
Therefore, the claims (i) and (ii) can be implied
immediately from the definition of q-module concavity or preordering concavity
and plugging in the value of $u$.
\end{proof}

\begin{theorem} \label{P:RepThm}
Assume $\mc{D}$ and $S_\mc{D}(f)$ are both convex and have nonempty interior.
\bdes
\item [(i)] If $f(x)$ and every $g_i(x)$ are q-module concave over $\mc{D}$,
then $S_\mc{D}(f) = \mc{L}_{qmod}^{\mc{D}}(f)$.

\item [(ii)] If $f(x)$ and every $g_i(x)$ are preordering concave over $\mc{D}$,
then $S_\mc{D}(f) = \mc{L}_{pre}^{\mc{D}}(f)$.
\edes
\end{theorem}
\begin{proof}
(i) Since $S_\mc{D}(f)$ is contained in the projection of $\mc{L}_{qmod}^{\mc{D}}(f)$,
we only need prove $S_\mc{D}(f) \supseteq \mc{L}_{qmod}^{\mc{D}}(f)$ .
For a contradiction, suppose there exists some
$\hat y \in \mc{L}_{qmod}^{\mc{D}}(f)$ such that $\hat x = \rho(\hat y) \notin S_\mc{D}(f)$.
By  the convexity of $S_\mc{D}(f)$, it holds
\[
S_\mc{D}(f) =  \bigcap_{\substack{\{a^Tx = b\} \mbox{ is a } \\ \mbox{ supporting hyperplane } }}  \quad
\Big\{x\in \re^n:\, a^Tx \geq b\Big\}.
\]
If $\hat x \notin S_\mc{D}(f)$, then there exists one hyperplane $\{a^Tx = b\}$ of $S_\mc{D}(f)$ such that
$a^T \hat x < b.$
%Consider optimization problem
%\begin{align*}
%b = \,\,\,\, \min_x &\quad a^Tx  \\
%\mbox{ subject to } & \quad f(x) \geq 0, \,\,\, g_1(x) \geq 0, \cdots, g_m(x)\geq 0.
%\end{align*}
By Lemma~\ref{lem:qplinrep},
we have representation
\be  \label{eq:Prep1}
a^Tx - b  =\lmd f(x)+ \sum_{i=0}^m g_i(x)\sig_i(x)
\ee
for some sos polynomials $\sig_i(x)$ such that
$\deg(g_i\sig_i) \leq 2 d_{qmod}^{(P)}$. Note that $d_{pre}^{(P)}=d$.
Write $\sig_i(x)$ as
\[
\sig_i(x) = \mathfrak{m}_{d-d_i}(x)^T W_i \mathfrak{m}_{d-d_i}(x),\, i=0,1,\ldots,m
\]
for some symmetric matrices $W_i \succeq 0$.
Then the identity \reff{eq:Prep1} becomes (noting \reff{P:mprod})
\[
a^Tx - b  =\lmd f(x)+  
\sum_{i=0}^m  (g_i(x)  \mathfrak{m}_{d-d_i}(x) \mathfrak{m}_{d-d_i}(x)^T )  \bullet W_i 
= \lmd f(x) + \sum_{i=0}^m 
\left( \sum_{\af\in\N^n: |\af| \leq 2d}  G_\af^{(i)} x^\af \right) \bullet W_i.
\]
In the above identity,
if we replace each $x^\af$ by $\hat y_\af$,
then get the contradiction
\[
a^T \hat x - b = f^T \hat y + \sum_{i=0}^m N_i(\hat y) \bullet W_i \geq 0.
\]
 
(ii) The proof is almost the same as for (i). The only difference is that
we have a new representation
\[
a^Tx - b  =\lmd f(x)+ \sum_{\nu\in\{0,1\}^m} g_1^{\nu_1}(x)\cdots g_m^{\nu_m}(x)\sig_\nu(x)
\]
for some sos polynomials $\sig_\nu(x)$ such that
$\deg(g_1^{\nu_1}\cdots g_m^{\nu_m}\sig_\nu) \leq 2d_{pre}^{(P)}=2d$.
Write $\sig_i(x)$ as
\[
\sig_\nu(x) = \mathfrak{m}_{d-d_\nu}(x)^T W_\nu \mathfrak{m}_{d-d_\nu}(x), \, W_\nu \succeq 0.
\]
%\begin{align*}
%b \quad = \qquad  \max  & \quad \gamma \\
% \mbox{ subject to } &  \quad
%a^Tx - \gamma = \lmd f(x)+ \sum_{\nu\in\{0,1\}^m} g_1^{\nu_1}(x)\cdots g_m^{\nu_m}(x)\sig_\nu(x),  \\
% & \qquad \qquad \lmd \geq 0,  \mbox{ every } \sig_\nu(x) \mbox{ is sos. }
%\end{align*}
%The dual of the above convex optimization problem is
%\begin{align*}
%  \min  \quad a^Tx  \qquad
%\mbox{ subject to }  \, \, x = \rho(y), \, y \in \mc{L}_{pre}^{\mc{D}}(f).
%\end{align*}
Then a similar contradiction argument can be applied prove the claim.
\end{proof}

\subsection{Some special cases}

Now we turn to some special cases about
q-module or preordering convexity/concavity or SDP representations.

\subsubsection{The q-module or preordering convexity certificate using Hessian}

The q-module or preordering convexity of $f(x)$
over the domain $\mc{D}$
can be verified by solving some semidefinite programming.
See \cite{ParMP, Las01} about the sos polynomials and semidefinite programming.
However, in some special cases like $\mc{D}=\re_+^n$,
a certificate for semidefiniteness of the Hessian $\nabla^2 f(x)$
can be applied to prove the q-module or preordering convexity of $f(x)$.

First, consider the case that $g_k(x)$ are concave over $\re^n$.  By concavity, it holds
\[
g_k(sx + (1-s)u) \geq s g_k(x) + (1-s) g_k(u), \quad \forall \,\, s \in [0,1], \, x, u \in \re^n.
\]
Now we assume the following certificate for the above criteria
\be  \label{eq:OmegaSOS}
\baray{c}
g_k(sx + (1-s)u) - s g_k(x) - (1-s) g_k(u) =  \\
\sig_0^{(k)}(x,u,s) + s\sig_1^{(k)}(x,u,s) + (1-s)\sig_2^{(k)}(x,u,s) + s(1-s) \sig_3^{(k)}(x,u,s)
\earay
\ee
where $\sig_0^{(k)}(x,u,s),\sig_1^{(k)}(x,u,s), \sig_2^{(k)}(x,u,s),\sig_3^{(k)}(x,u,s)$
are sos polynomials in $(x,u,s)$.
Note that the identity \reff{eq:OmegaSOS} is always true
when $\mc{D}$ is a polyhedra, i.e., every $g_k(x)$ has degree one.

\begin{theorem} \label{thm:HesCrtf}
Suppose for every $1\leq k \leq m$, the identity \reff{eq:OmegaSOS} holds.
If $\nabla^2f(x)$ belongs to the quadratic module (resp. preordering)
generated by polynomials $g_1(x), \cdots, g_m(x)$, i.e.,
\[
\nabla^2 f(x) =  \sum_{i=0}^m  g_i(x) H_i(x) \,\,
\left(resp. \quad  \nabla^2 f(x) =  \sum_{\nu\in \{0,1\}^m }  
g_1^{\nu_1}(x) \cdots g_m^{\nu_m}(x) H_\nu(x)  \right) 
\]
for some sos matrices $H_i(x)$ (resp. $H_\nu(x)$),
then $f(x)$ is q-module convex (resp. preordering convex)
over the domain $\mc{D}$.
\end{theorem}
\begin{proof}
First suppose $\nabla^2f(x)$ belongs to the quadratic module
generated by $g_1(x), \cdots, g_m(x)$, i.e.,
$ \nabla^2 f(x) = \sum_{i=0}^m g_i(x) H_i(x)$ (recall $g_0(x)=1$)
for some sos matrices $H_i(x)$.
Then we have
\begin{align*}
& \, f(x) - f(u) - \nabla f(u)^T(x-u) \\
 = &  (x-u)^T \left( \int_0^1 \int_0^t  \nabla^2 f( u + s(x-u)) \,\, ds \, dt  \right) (x-u) \\
 = &  (x-u)^T \left( \sum_{k=1}^m  \int_0^1 \int_0^t g_k( sx + (1-s)u ) H_i( sx +(1-s)u )  \,\, ds \, dt \right) (x-u).
\end{align*}
By identity \reff{eq:OmegaSOS}, we have
\begin{align*}
& \, \quad \int_0^1 \int_0^t g_k( sx + (1-s)u ) H_i( sx +(1-s)u )  \,\, ds \, dt  \\
= & \int_0^1 \int_0^t  \Big( \sig_0^{(k)}(x,u,s)+ \sig_1^{(k)}(x,u,s) + 
(1-s)\sig_2^{(k)}(x,u,s) + s(1-s) \sig_3^{(k)}(x,u,s) \Big)
H_i( sx +(1-s)u )  \,\, ds \, dt \\
& + g_k(x) \int_0^1 \int_0^t s H_i( sx +(1-s)u )  \,\, ds \, dt
+ g_k(u) \int_0^1 \int_0^t (1-s) H_i( sx +(1-s)u )  \,\, ds \, dt \\
= & H_0^{(k)}(x,u) +  g_k(x)H_1^{(k)}(x,u) + g_k(u)H_2^{(k)}(x,u)
\end{align*}
for some sos matrices $H_0^{(k)}(x,u), H_1^{(k)}(x,u),H_2^{(k)}(x,u)$
(see Lemma~3.1 in \cite{HN1}).
So $f(x)$ is q-module convex over $\mc{D}$.

\medskip

Second,  when $\nabla^2f(x)$ belongs to the preordering
generated by $g_1(x), \cdots, g_m(x)$, i.e.,
\[
\nabla^2 f(x) =  \sum_{\nu\in \{0,1\}^m} g_1^{\nu_1}(x)\cdots g_m^{\nu_m}(x) H_\nu(x)
\]
for some sos matrices $H_\nu(x)$,
a similar argument as above shows $f(x)$ is preordering convex over $\mc{D}$.
\end{proof}

Second, consider the special case that $\mc{D} = \re_+^n$
and the polynomial $f(x)$ is cubic.

\begin{theorem} \label{thm:cubR+3}
Let $\re_+^n=\{x:x_1\geq0,\cdots, x_n \geq 0\}$ be the domain.
If $f(x)$ is a cubic polynomial concave over $\re_+^n$,
then $S_{\re_+^n}(f) = \mc{L}_{qmod}^{\re_+^n}(f)$.
\end{theorem}
\begin{proof}
By Theorem~\ref{P:RepThm}, it suffices to
prove $f(x)$ is q-module concave over $\re_+^n$.
Since $f(x)$ is cubic, we have
\[
- \nabla^2 f(x) = A_0 +x_1A_1 + \cdots + x_nA_n
\]
for some symmetric matrices $A_i$.
When $f(x)$ is concave over $\re_+^n$,
$- \nabla^2 f(x) \succeq 0$ for all $x\geq 0$.
Hence all $A_i$ must be positive semidefinite.
This means that $-\nabla^2 f(x)$ belongs to the quadratic module
generated by $x_1,\cdots,x_n$.
By Theorem~\ref{thm:HesCrtf},
$f(x)$ is q-module concave over $\re_+^n$.
%So we can write $A_i = B_i^TB_i$.
%Therefore, it holds
%\begin{align*}
%& \, f(u) + \nabla f(u)^T(x-u) - f(x) \\
% = &  (x-u)^T \left( \int_0^1 \int_0^t A_0 + \sum_{i=1}^n (u_i+s(x_i-u_i)) A_n  \right) (x-u) \\
% = & \half \|B_0(x-u)\|^2 + \sum_{i=1}^n  \left( \frac 13 u_i \|B_i (x-u)\|^2
% + \frac 16 x_i \|B_i (x-u)\|^2 \right).
%\end{align*}
%So $f(x)$ is q-module concave over $\re_+^n$.
\end{proof}

\begin{example}
The convex set $S_{\re_+^2}(f)$
with $f(x)=1-(x_1^3+x_1^2x_2+x_1x_2^2+x_2^3)$
equals $\mc{L}_{qmod}^{\re_+^2}(f)$:
\[
\left\{(x_1,x_2): \,\exists \, y_{ij}, s.t. \, \,
\baray{c}
1 \geq y_{30} + y_{21} + y_{12} + y_{03}   \\
\bbm
1  &  x_1  & x_2  & y_{20} & y_{11} & y_{02} \\
x_1  &  y_{20}  & y_{11}  & y_{30} & y_{21} & y_{12} \\
x_2 &  y_{11}  & y_{02} & y_{21} & y_{12} & y_{03} \\
y_{20}  &  y_{30}  & y_{21}  & y_{40} & y_{31} & y_{22} \\
y_{11} &  y_{21}  & y_{12} & y_{31} & y_{22} & y_{13} \\
y_{02} &  y_{12}  & y_{03} & y_{22} & y_{13} & y_{04} \\
\ebm \succeq 0,
\baray{c}
\bbm
x_1  &  y_{20}  & y_{11}  \\
y_{20}  &  y_{30}  & y_{21}  \\
y_{11} &  y_{21}  & y_{12}
\ebm \succeq 0 \\
\bbm
x_2  &  y_{11}  &y_{02}  \\
y_{11}  &  y_{21}  & y_{12}  \\
y_{02} &  y_{12}  & y_{03}
\ebm \succeq 0
\earay
\earay
\right\}.
\]
This is because $f(x)$ is q-mod concave over $\re_+^2$
(see Example~\ref{exm:2dot1}).
\end{example}

\medskip

Third, consider the special case of univariate polynomials.
When $\mc{D}= \re$,
a univariate polynomial is convex if and only if it is sos-convex,
which holds if and only if it is first order sos convex.
When $\mc{D}= I$ is an interval,
we will see that a univariate polynomial is convex over $I$
if and only if it is q-module convex over $I$.

%
%
%For example, the univariate polynomial $f(x)=x^d$ ($d>2$ is odd) is not globally convex,
%but is convex over $R_+$.
%We are interested to know whether it is also q-module convex with respect to $\re_+$.
%When $d=3$, the q-module convexity is obtained due to the identity
%\[
%x^3-u^3-3u^2(x-u) = x (x-u)^2 + 2 u (u-x)^2.
%\]
%When $d=5$, the q-module convexity can be shown by observing the identity
%\[
%x^5-u^5-5u^4(x-u) = x (x(x-u))^2 + 2 u (x(x-u))^2 + 3x(u(x-u))^2+4(u(x-u))^2.
%\]
%Actually, $x^d$ is always q-module convex over $\re_+$
%for any positive integer $d$.

\begin{pro}  \label{pro:upolyCvx}
Let $f(x)$ be a univariate polynomial,
and $I$ be an interval like $[a,b]$, $(-\infty, b]$ or $[a,\infty)$.
Then $f(x)$ is convex over $I$
if and only if it is q-module convex over $I$.
\end{pro}
\begin{proof}
First suppose $I=[a,b]$ is finite. $f(x)$ is convex over $[a,b]$ if and only if
$f^{\prm\prm}(x)\geq 0$ for all $x\in [a,b]$, which is true if and only if
\[
f^{\prm\prm}(x) = \sig_0(x) + (x-a)\sig_1(x) + (b-x) \sig_2(x)
\]
for some sos polynomials $\sig_0, \sig_1, \sig_2$ with degrees at most
$2\lceil deg(f)/2 \rceil$ (see Powers and Reznick \cite{PR00}).
In other words, $f(x)$ is convex over $[a,b]$ if and only if its Hessian
belongs to the quadratic module generated by polynomials $x-a, b-x$.
Then the conclusion can be implied by Theorem~\ref{thm:HesCrtf}.

The proof is similar for the case $(-\infty, b]$ or $[a,\infty)$.
\end{proof}

\subsubsection{Epigraph of polynomial functions}

For a given convex domain $\mc{D} \subseteq \re^n$,
$f(x)$ is convex over $\mc{D}$
if and only if its epigraph
\[
\mathbf{epi}(f) := \{(x,t)\in \mc{D} \times \re: \,
f(x) \leq t \}
\]
is convex.
Note that $\mathbf{epi}(f)$ is defined by the inequality
$ t-f(x) \geq 0 $.
If we consider $t-f(x)$ as a polynomial in $x$
with coefficients in $t$,
then $\mc{L}_{qmod}^{\mc{D}}(t-f)$
and $\mc{L}_{pre}^{\mc{D}}(t-f)$ are both linear in $(x,t)$.
Therefore, if $f(x)$ is q-module (resp. preordering) convex over $\mc{D}$,
$\mc{L}_{qmod}^{\mc{D}}(t-f)$ (resp. $\mc{L}_{pre}^{\mc{D}}(t-f)$)
presents an SDP representation for $\mathbf{epi}(f)$.

By Proposition~\ref{pro:upolyCvx},
when $f(x)$ is a univariate polynomial convex over
an interval $I$,
we know its epigraph $ \mbox{epi} (f) $
is SDP representable
and $\mc{L}_{qmod}^{I}(t-f)$ is one SDP representation.

\medskip

\section{Convex sets defined by rational functions} 
\label{sec:RatCvx}

In this section, we discuss the SDP representation of convex set
$S_\mc{D}(f)$ when $f(x)$ is a rational function
while the domain $\mc{D}=\{x\in\re^n:g_1(x)\geq 0, \cdots, g_m(x) \geq 0 \}$
is still defined by polynomials $g_i$.
Let $f(x)$ be a rational function of the form
\[
f(x)= \frac {1}{f_{den}(x)} \sum_{\af \in \N^n:\,  |\af| \leq 2d}   f_\af x^\af.
\]
Here $f_{den}(x)$ is the denominator of $f(x)$.
We assume that $f(x)$ is concave over the domain $\mc{D}$.
So $f(x)$ can not have poles in the interior $int(\mc{D})$ of $\mc{D}$.
Without loss of generality, assume $f_{den}(x)$ is positive over $int(\mc{D})$.
Note that $f(x)$ is not defined
on the boundary $\pt \mc{D}$ where $f_{den}(x)$ vanishes.
If this happens, we think of $ S_\mc{D}(f)$
as the closure of $\{x\in int(\mc{D}):\, f(x) \geq 0 \}$.

\subsection{The q-module or preordering convexity of rational functions}

We now introduce some types of definitions about convexity/concavity for rational functions.
Let $p(x),q(x)$ be two given polynomials
which are positive in $int(\mc{D})$.
We say $f(x)$ is {\it q-module convex over $\mc{D}$ with respect to $(p,q)$}
if the identity
\be \label{eq:qmodR}
p(x)q(u)\cdot R_f(x,u) = \sum_{i=0}^m g_i(x) \left( \sum_{j=0}^m g_j(u) \sig_{ij}(x,u)  \right)
\ee
holds for some sos polynomials $\sig_{ij}(x,u)$.
Then define $f(x)$ to be {\it q-module concave over $\mc{D}$ with respect to $(p,q)$}
if $-f(x)$ is  q-module convex over $\mc{D}$ with respect to $(p,q)$.
We say $f(x)$ is {\it preordering convex over $\mc{D}$ with respect to $(p,q)$}
if the identity
\be\label{eq:preR}
p(x)q(u) \cdot R_f(x,u) = \sum_{\nu\in \{0,1\}^m } g_1^{\nu_1}(x)\cdots g_m^{\nu_m}(x)
\left(  \sum_{\mu\in \{0,1\}^m } g_1^{\mu_1}(u)\cdots g_m^{\mu_m}(u)  \sig_{\nu,\mu}(x,u)  \right)
\ee
holds for some sos polynomials $\sig_{\nu,\mu}(x,u)$.
Similarly, $f(x)$ is called {\it preordering concave over $\mc{D}$ with respect to $(p,q)$}
if $-f(x)$ is preordering convex over $\mc{D}$ with respect to $(p,q)$.
We point out that the definition of 
q-module or preordering convexity/concavity over $\mc{D}$
for rational functions assumes a certain set of defining polynomials 
$g_i(x)$ for $\mc{D}$ is clear in the context.

In identities \reff{eq:qmodR} or \reff{eq:preR},
there is no information on how to find polynomials $p,q$.
However, since $R_f(x,u)$ has denominator $f_{den}(x)f_{den}^2(u)$,
a possible choice for $(p,q)$ is
\be \label{eq:denpq}
p(x) = f_{den}(x), \quad q(u) = f_{den}^2(u).
\ee
If the choice $(p,q)$ in \reff {eq:denpq} makes
the identity \reff{eq:qmodR} (resp. \reff{eq:preR}) holds,
we say $f(x)$ is q-module (resp. preordering)
convex over $\mc{D}$ with respect to $p(x)$,
or just simply say $f(x)$ is q-module (resp. preordering)
convex over $\mc{D}$
if the denominator $f_{den}(x)$ is clear in the context.

In the special case $\mc{D} = \re^n$,
the definitions of q-module and preordering convexity over $\mc{D}$
coincide with each other,
and then is called {\it first order sos convexity}
when $(p,q)$ is given by \reff {eq:denpq},
as consistent with the definition of first order sos convexity
in Section~\ref{sec:poly}.
First order sos concavity is defined similarly.

\begin{example}
(i) The rational function $\frac{x_2^2}{x_1}$ is convex over
the domain $\re_+\times \re$.  It is also q-module convex over $\re_+\times \re$
with respect to the denominator $x_1$, which is due to that
\[
\frac{x_2^2}{x_1}-\frac{u_2^2}{u_1}
-\left(-\frac{u_2^2}{u_1^2}(x_1-u_1)+\frac{2u_2}{u_1}(x_2-u_2) \right)
= \frac{1}{x_1u_1^2}(x_1u_2-x_2u_1)^2.
\]
%Obviously it is also preordering convex over $\re_+\times \re$.
%A more general case is that
%\begin{align*}
%\frac{x_2^t}{x_1^s}-\frac{u_2^t}{u_1^s}
%-\left(-s\frac{u_2^t}{u_1^{s+1}}(x_1-u_1)+\frac{tu_2^{t-1}}{u_1^s}(x_2-u_2) \right) = \hspace{1cm} \\
%\frac{1}{x_1^su_1^{s+1}}
%\left(x_2^tu_1^{s+1} + (t-s-1)x_1^su_1u_2^t+sx_1^{s+1}u_2^t-tx_1^sx_2u_1u_2^{t-1} \right).
%\end{align*}
\\
(ii) The rational function $f(x)=\frac{x_1^4+x_1^2x_2^2+x_2^4}{x_1^2+x_2^2}$ is convex over
the domain $\re^2$. It can be verified that
\[
(f(x) - f(u) - \nabla f(u)(x-u)) =
\frac{f_1^2 + f_2^2 + \half (f_3^2 + f_4^2 + f_5^2 + f_6^2) + f_7^2 + f_8^2 + f_9^2}
{(x_1^2+x_2^2)(u_1^2+u_2^2)^2}
\]
where the polynomials $f_i$ are given as below
\begin{align*}
\baray{ll}
f_1 = -u_1 u_2 x_2^2-u_1 u_2 x_1^2+u_1 u_2^2 x_2+u_1^2 u_2 x_1, &
f_6 = -u_2^2 x_2^2+u_2^3 x_2-u_1^2 x_1^2+u_1^3 x_1,\\
f_2 = -u_1 u_2 x_2^2+u_1 u_2 x_1^2+u_1 u_2^2 x_2-u_1^2 u_2 x_1, &
f_7 = -2 u_1 u_2 x_1 x_2 + u_1 u_2^2 x_1 +  u_1^2 u_2 x_2,\\
f_3 = -u_2^2 x_1 x_2+u_2^3 x_1-u_1^2 x_1 x_2+u_1^3 x_2, &
f_8 = u_2^2 x_1^2 - u_1^2 x_2^2, \\
f_4 =  u_2^2 x_1 x_2-u_2^3 x_1-u_1^2 x_1 x_2+u_1^3 x_2, &
f_9 = -u_1 u_2^2 x_1+u_1^2 u_2 x_2. \\
f_5 =  u_2^2 x_2^2-u_2^3 x_2-u_1^2 x_1^2+u_1^3 x_1,\\
\earay
%f_6 &= -u_2^2 x_2^2+u_2^3 x_2-u_1^2 x_1^2+u_1^3 x_1,\\
%f_7 &= -2 u_1 u_2 x_1 x_2 + u_1 u_2^2 x_1 +  u_1^2 u_2 x_2,\\
%f_8 &= u_2^2 x_1^2 - u_1^2 x_2^2, \\
%f_9 &= -u_1 u_2^2 x_1+u_1^2 u_2 x_2.
\end{align*}
So the $f(x)$ given above is first order sos convex. 
\qed
\end{example}

\medskip
Obviously, the q-module convexity implies preordering convexity,
which then implies the convexity,
but the converse might not be true.
For instance, $\frac{1}{x_1x_2}$ is convex over $\re_+^2$, but it is
neither q-module nor preordering convex over $\re_+^2$
with respect to the denominator $x_1x_2$. 
Note that for $f(x) = \frac{1}{x_1x_2}$ it holds
\[
u_1^2u_2^2x_1x_2 R_f(x,u) =
u_1^2u_2^2+x_1^2x_2u_2+x_1x_2^2u_1-3x_1x_2u_1u_2.
\]
There are no sos polynomials $\sig_{\nu,\mu}(x,u)$ such that
\[
u_1^2u_2^2x_1x_2 R_f(x,u) =  \sum_{\nu\in \{0,1\}^2, \mu\in \{0,1\}^2} 
x_1^{\nu_1}x_2^{\nu_2} u_1^{\mu_1}u_2^{\mu_2} \sig_{\nu,\mu}(x,u).
%\sig_0 + x_1x_2 \sig_1 + u_1u_2 \sig_2
%+x_1u_1 \sig_3 + x_1u_2 \sig_4 + x_2u_1 \sig_5 + x_2u_2 \sig_6.
\]
Otherwise, if they exist, we replace
$(x_1,x_2,u_1,u_2)$ by $(x_1^2,x_2^2,1,1)$ and then get the dehomogenized Motzkin's polynomial
$1+x_1^4x_2^2+x_1^2x_2^4-3x_1^2x_2^2$
is sos, which is impossible (see Reznick \cite{Rez00}).

\begin{pro}
Let $CQ_{p,q}(\mc{D})$ (resp. $CP_{p,q}(\mc{D})$) be the set of all q-module
(resp. preordering) convex rational functions
over $\mc{D}$ with respect to $(p,q)$.
Then they have the properties:
\bit
\item [(i)] Both $CQ_{p,q}(\mc{D})$ and $CP_{p,q}(\mc{D})$ are convex cones.

\item [(ii)] If $f(x) \in CQ_{p,q}(\mc{D})$ (resp. $f(x) \in CP_{p,q}(\mc{D})$),
then $f(Az+b) \in CQ_{\tilde p, \tilde q}(\tilde{\mc{D}})$
(resp. $f(Az+b) \in CP_{\tilde p, \tilde q}(\tilde{\mc{D}})$),
where $\tilde{\mc{D}} = \{z : Az+b\in \mc{D} \}$
and $\tilde p(z) = p(Az+b), \tilde q(z) = q(Az+b)$.
That is, the q-module convexity or preordering convexity is preserved
under linear transformations.
\eit
\end{pro}
\begin{proof}
The item (i) can be verified explicitly,
and item (ii) can be obtained by substituting $Az+b$ for $x$
and noting the chain rule of derivatives.
\end{proof}

\subsection{SDP representations}

Now we turn to the construction of SDP representations for $S_\mc{D}(f)$.
Recall $g_0(x) \equiv 1$.
Throughout this subsection, we assume the polynomials $f(x)$
and $g_i(x)$ are either all q-module concave or all preordering concave over $\mc{D}$
with respect to $(p,q)$.
Thus for any $h(x)$ from $\{f(x), g_1(x), \cdots, g_m(x)\}$, we have either 
\[
-p(x)q(u)R_h(x,u) = \sum_{i=0}^m g_i(x) \left( \sum_{j=0}^m g_j(u) \sig_{ij}^{h}(x,u)  \right)
\]
for some sos polynomials $\sig_{ij}^{h}(x,u)$, or 
\[
-p(x)q(u)R_h(x,u) = \sum_{\nu\in \{0,1\}^m } g_1^{\nu_1}(x)\cdots g_m^{\nu_m}(x)
\left(  \sum_{\mu\in \{0,1\}^m } g_1^{\mu_1}(u)\cdots g_m^{\mu_m}(u)  \sig_{\nu,\mu}^{h}(x,u)  \right)
\]
for some sos polynomials $\sig_{\nu,\mu}^{h}(x,u)$. Let
\be \label{R:dgfrm}
\baray{c}
d_{qmod}^{(R)} = \underset{h \in \{f, g_1, \cdots, g_m\}}{\max} \quad
\underset{0\leq i,j \leq m}{\max} \quad
\lceil \half \deg_x(g_i \sig_{ij}^{h}) \rceil,\\
d_{pre}^{(R)} = \underset{h \in \{f, g_1, \cdots, g_m\}}{\max} \quad
\underset{\nu\in \{0,1\}^m, \mu\in \{0,1\}^m }{\max} \quad
\lceil \half \deg_x(g_1^{\nu_1}\cdots g_m^{\nu_m} \sig_{\nu,\mu}^{h}) \rceil.
\earay
\ee
Then set $d=d_{qmod}^{(R)}$ (resp. $d=d_{pre}^{(R)}$)
when $f(x)$ and $g_i(x)$ are all q-module (resp. preordering) concave over $\mc{D}$
with respect to $(p,q)$.

Define matrices $P_\af^{(i)}, P_\af^{(\nu)}, Q_\af^{(i)}, Q_\af^{(\nu)}$ such that
\be \label{R:mprod}
\baray{rcl}
\frac{g_i(x)}{p(x)} \mathfrak{m}_{d-d_i}(x)\mathfrak{m}_{d-d_i}(x)^T &= &
\underset{\af \in \N^n: |\af| + |LE(p)| \leq 2d }{\sum} Q_\af^{(i)} x^\af
+ \underset{\bt \in \N^n: \bt < LE(p)}{\sum} P_\af^{(i)} \frac{x^\bt}{p(x)}, \quad 0\leq i\leq m, \\
\frac{g_1^{\nu_1}\cdots g_m^{\nu_m}}{p(x)}
\mathfrak{m}_{d-d_\nu}(x)\mathfrak{m}_{d-d_\nu}(x)^T &= &
\underset{\af \in \N^n: |\af| + |LE(p)| \leq 2d }{\sum} Q_\af^{(\nu)} x^\af
+ \underset{\bt \in \N^n: \bt < LE(p)}{\sum} P_\af^{(\nu)} \frac{x^\bt}{p(x)}, \quad \nu\in \{0,1\}^m.
\earay
\ee
Here $LE(p)$ denotes the exponent of
the leading monomial of $p(x)$
under the lexicographical ordering
($x_1 > x_2  > \cdots > x_n $).
Note that the union
\[
\Big\{x^\af:  \af \in \N^n, |\af| + |LE(p)| \leq 2d \Big\} \cup
\Big\{\frac{x^\bt}{p(x)}:  \bt \in \N^n, \bt < LE(p)\Big\}
\]
is a set of polynomials and rational functions that are linearly independent.

Let $y$ be a vector indexed by $\af \in \N^n$ such that $|\af| + |LE(p)| \leq 2d$,
and $z$ be a vector indexed by $\bt\in \N^n$ such that $\bt < LE(p)$.
Then define
\be  \label{lmi:Qyz}
\baray{l}
Q_i(y,z) =
\underset{\af \in \N^n: |\af| + |LE(p)| \leq 2d }{\sum} Q_\af^{(i)} y_\af
+ \underset{\bt \in \N^n: \bt < LE(p)}{\sum} P_\af^{(i)}  z_\bt, \, 0\leq i\leq m,  \\
Q_\nu(y,z)=
\underset{\af \in \N^n: |\af| + |LE(p)| \leq 2d }{\sum} Q_\af^{(\nu)} y_\af
+\underset{\bt \in \N^n: \bt < LE(p)}{\sum} P_\af^{(\nu)} z_\bt, \nu \in \{0,1\}^m.
\earay
\ee
Suppose the rational function $f(x)$ is given in the form
\[
f(x) =
\sum_{\af \in \N^n: |\af| + |LE(p)| \leq 2d } f_\af^{(1)} x^\af
+ \sum_{\bt \in \N^n: \bt < LE(p)} f_\bt^{(2)} \frac{x^\bt}{p(x)},
\]
then define vectors $f^{(1)}, f^{(2)}$ such that
\[
(f^{(1)})^Ty + (f^{(2)})^Tz =
\sum_{\af \in \N^n: |\af| + |LE(p)| \leq 2d } f_\af^{(1)} y_\af
+ \sum_{\bt \in \N^n: \bt < LE(p)} f_\bt^{(2)} z_\bt.
\]
Define two SDP representable sets
\begin{align}
\mc{R}_{qmod}^{\mc{D}}(f) & =\left\{(y,z): \, y_{0} = 1, \,
\bbm f^{(1)} \\ f^{(2)} \ebm^T \bbm y \\ z \ebm \geq 0,\,
Q_i(y,z) \succeq 0, \, \forall \, 0\leq i \leq m \right\},  \label{LRqmod}\\
\mc{R}_{pre}^{\mc{D}}(f) & =\left\{(y,z): \, y_{0} = 1, \,
\bbm f^{(1)} \\ f^{(2)} \ebm^T \bbm y \\ z \ebm  \geq 0,\,
Q_\nu(y,z) \succeq 0, \, \forall \, \nu \in \{0,1\}^m \right\}. \label{LRpre}
\end{align}
%Then the set $S_{\mc{D}}(f)$ is contained in the projection of
%both $\mc{R}_{qmod}^{\mc{D}}(f)$ and $\mc{R}_{pre}^{\mc{D}}(f)$
%under the projection map:
%\[
%\rho(y,z) = (y_{1 0 \ldots 0}, \,\,\, y_{0 1 0 \ldots 0}, \,\,\, \ldots, \,\,\, y_{0 0 \ldots 0 1}).
%\]
We say $S_{\mc{D}}(f)$ equals $\mc{R}_{qmod}^{\mc{D}}(f)$ (resp. $\mc{R}_{pre}^{\mc{D}}(f)$)
if $S_{\mc{D}}(f)$ equals the image of $\mc{R}_{qmod}^{\mc{D}}(f)$ (resp. $\mc{R}_{pre}^{\mc{D}}(f)$)
under the projection
$
\rho(y,z) = (y_{1 0 \ldots 0}, \,\,\, y_{0 1 0 \ldots 0}, \,\,\, \ldots, \,\,\, y_{0 0 \ldots 0 1}).
$
Similarly, we say $x\in \mc{R}_{qmod}^{\mc{D}}(f)$ (resp. $x\in \mc{R}_{pre}^{\mc{D}}(f)$)
if there exists $(y,z) \in \mc{R}_{qmod}^{\mc{D}}(f)$ (resp. $(y,z)\in \mc{R}_{pre}^{\mc{D}}(f)$)
such that $x = \rho(y,z)$.

\begin{lemma}  \label{lem:RATlinrep}
Assume $\mc{D}$ and $S_\mc{D}(f)$ are both convex and have nonempty interior.
Let $\{x\in\re^n: a^Tx = b\} $ be a supporting hyperplane of $S_\mc{D}(f)$ such that
$a^Tx \geq b$ for all $x \in S_\mc{D}(f)$
and $a^Tu = b$ for some point $u\in S_\mc{D}(f) $ such that 
$q(u)>0$, and either $f(u)>0$ or $f_{den}(u)>0$.
\bdes
\item [(i)] If $f(x)$ and every $g_i(x)$ are q-module concave over $\mc{D}$
with respect to $(p,q)$, then
\[
p(x) \cdot (a^Tx - b - \lmd f(x)) =  \sum_{i=0}^m  g_i(x) \sig_i(x)
\]
for some scalar $\lmd \geq 0$ and sos polynomials $\sig_i(x)$ 
such that $\deg(g_i\sig_i) \leq 2d_{qmod}^{(R)}$.

\item [(ii)] If $f(x)$ and every $g_i(x)$ are preordering concave over $\mc{D}$
with respect to $(p,q)$, then
\[
p(x) \cdot (a^Tx - b - \lmd f(x)) =  \sum_{\nu \in \{0,1\}^m }  
g_1^{\nu_1}(x) \cdots g_m^{\nu_m}(x) \sig_\nu(x)
\]
for some scalar $\lmd \geq 0$ and sos polynomials $\sig_i(x)$ 
such that $\deg(g_1^{\nu_1} \cdots g_m^{\nu_m} \sig_\nu) \leq 2d_{pre}^{(R)}$.
\edes
\end{lemma}
\begin{proof}
Since $S_\mc{D}(f)$ has nonempty interior, the first order optimality condition holds
at $u$ for convex optimization problem ($u$ is one minimizer)
\[
\min_{x}  \quad a^Tx \quad
\mbox{ subject to } \quad  f(x)\geq 0, \,  g_i(x)\geq 0,  \, i=1, \ldots, m.
\]
If $f(u)>0$, the constraint $f(x)\geq 0$ is inactive.
If $f_{den}(u)>0$, $f(x)$ is differential at $u$.
Hence, in either case, there exist Lagrange multipliers
$\lmd\geq 0, \lmd_1 \geq 0, \ldots, \lmd_m \geq 0$ such that
\begin{align*}
a =  \lmd \nabla f(u) + \sum_{i=1}^m \lmd_i \nabla g_i(u), \quad
\lmd f(u) = \lmd_1 g_1(u) = \cdots =\lmd_m g_m(u) = 0.
\end{align*}
Hence we get the identity
\begin{align*}
a^Tx - b - \lmd f(x)-\sum_{i=1}^m \lmd_i g_i(x) & = \\
\lmd \Big ( f(u) + \nabla f(u)^T(x-u)  - f(x) \Big) +
\sum_{i=1}^m \lmd_i \Big ( g_i(u) & + \nabla g_i(u)^T(x-u)  - g_i(x) \Big) .
\end{align*}
Therefore, the claims (i) and (ii) can be implied
immediately by the definition of q-module concavity or preordering concavity of 
$f$ and $g_i$, and plugging the value of $u$.
\end{proof}

\begin{theorem}  \label{thm:RatRep}
Assume $\mc{D}$ and $S_\mc{D}(f)$ are both convex and have nonempty interior.
Let $(p,q)$ be given in \reff{eq:qmodR} or \reff{eq:preR}.
Suppose $\dim(\mc{Z}(f) \cap \mc{Z}(f_{den}) \cap \pt S_\mc{D}(f)) < n-1$ and
$\dim(\mc{Z}(q) \cap \pt S_\mc{D}(f)) < n-1$.
\bdes
\item [(i)] If $f(x)$ and every $g_i(x)$ are q-module concave over $\mc{D}$ with respect to $(p,q)$, then
$S_\mc{D}(f) = \mc{R}_{qmod}^{\mc{D}}(f)$.

\item [(ii)] If $f(x)$ and every $g_i(x)$ are preordering concave over $\mc{D}$ with respect to $(p,q)$, then
$S_\mc{D}(f) = \mc{R}_{pre}^{\mc{D}}(f)$.

\edes
\end{theorem}
\begin{proof}
(i) Since $S_\mc{D}(f)$ is contained in the projection of $\mc{R}_{qmod}^{\mc{D}}(f)$,
we only need prove $S_\mc{D}(f) \supseteq \mc{R}_{qmod}^{\mc{D}}(f)$.
For a contradiction, suppose there exists some
$(\hat y, \hat z) \in \mc{R}_{qmod}^{\mc{D}}(f)$ such that $\hat x = \rho(\hat y) \notin S_\mc{D}(f)$.
By the convexity of $S_\mc{D}(f)$, it holds
\[
S_\mc{D}(f) =  \bigcap_{\substack{\{a^Tx = b\} \mbox{ is a } \\ \mbox{ supporting hyperplane } }}  \quad
\Big\{x\in \re^n:\, a^Tx \geq b\Big\}.
\]
If $\hat x \notin S_\mc{D}(f)$, then there exists one supporting hyperplane
$\{a^Tx = b\}$ of $S_\mc{D}(f)$ with tangent point $u \in \pt S_\mc{D}(f)$ such that
$a^T \hat x < b.$
Since $\dim(\mc{Z}(f) \cap \mc{Z}(p) \cap \pt S_\mc{D}(f)) < n-1$
and $\dim(\mc{Z}(q) \cap \pt S_\mc{D}) < n-1$,
by continuity, we can choose
$\{a^Tx = b\}$ such that either $f(u)>0$ or $p(u) > 0$, 
and $q(u)>0$.
%Consider optimization problem
%\begin{align*}
%b = \,\,\,\, \min_x &\quad a^Tx  \\
%\mbox{ subject to } & \quad f(x) \geq 0, \,\,\, g_1(x) \geq 0, \cdots, g_m(x)\geq 0.
%\end{align*}
By Lemma~\ref{lem:RATlinrep}, we have
\be \label{eq:Rlrep2}
a^Tx - b  =\lmd f(x)+ \sum_{i=0}^m \frac{g_i(x)}{p(x)} \sig_i(x)
\ee
for some sos polynomials $\sig_i(x)$ such that $\deg(g_i \sig_i) \leq 2 d_{qmod}^{(R)}$.
Note that $d = d_{qmod}^{(R)}$. Then write $\sig_i(x)$ as
\[
\sig_i(x) = \mathfrak{m}_{d-d_i}(x)^T W_i \mathfrak{m}_{d-d_i}(x),\, i=0,1,\ldots,m
\]
for some symmetric matrices $W_i \succeq 0$.
%Now, if we choose vectors $\bar y$ and $\bar z$ as
%\[
%\bar y_\af =  \hat x^\af, \af: |\af| \leq 2d - |LE(p)|, \qquad
%\bar z_\bt = \frac{1}{p(\hat x)}  \hat x^\bt, \bt: \bt < LE(p),
%\]
Then identity \reff{eq:Rlrep2} becomes (noting \reff{R:mprod})
\begin{align*}
a^Tx - b  &=\lmd f(x)+ \sum_{i=0}^m 
\left(\frac{g_i(x)}{p(x)} \mathfrak{m}_{d-d_i}(x)\mathfrak{m}_{d-d_i}(x)^T\right) \bullet W_i \\
&= \lmd f(x) + 
\sum_{i=0}^m  \left(
\underset{\af \in \N^n: |\af| + |LE(p)| \leq 2d }{\sum}  Q_\af^{(i)} x^\af
+ \underset{\bt \in \N^n: \bt < LE(p)}{\sum} P_\af^{(i)}  \frac{x^\bt}{p(x)} \right) \bullet W_i.
\end{align*}
In the above identity, 
if we replace each $x^\af$ by $\hat y_\af$
and $\frac{x^\bt}{p(x)}$ by $\hat z_\bt$,
then we get the contradiction
\[
a^T \hat x - b = (f^{(1)})^T \hat y + (f^{(2)})^T \hat z  + 
\sum_{i=0}^m Q_i(\hat y, \hat z) \bullet W_i \geq 0.
\]
%Thus it must hold
%\begin{align*}
%b \quad = \qquad  \max  & \quad \gamma \\
% \mbox{ subject to } &  \quad
%a^Tx - \gamma = \lmd f(x)+ \sum_{i=0}^m \frac{g_i(x)}{p(x)}\sig_i(x),  \\
% & \qquad \qquad \lmd \geq 0,  \mbox{ every } \sig_i(x) \mbox{ is sos. }
%\end{align*}
%However, the dual of the above convex optimization problem is
%\begin{align*}
%  \min  \quad a^Tx  \qquad
%\mbox{ subject to }  \, \, x = \rho((y,w)), \, y \in \mc{L}_{qmod}^{\mc{D}}(f).
%\end{align*}
%By weak duality, it must have
%\[
%b \leq  a^Tx , \,\, \forall \, x = \rho((y,w)),  (y,w) \in \mc{R}_{qmod}^{\mc{D}}(f).
%\]

(ii) The proof is almost the same as for (i). The only difference is that
\[
a^Tx - \gamma = \lmd f(x)+ \sum_{\nu\in\{0,1\}^m} \frac{g_1^{\nu_1}\cdots g_m^{\nu_m}}{p(x)}\sig_\nu(x)
\]
for some sos polynomials $\sig_\nu(x)$ such that $\deg(g_i \sig_i) \leq 2 d_{pre}^{(R)}$.
Note that $d=d_{pre}^{(R)}$.
%\begin{align*}
%b \quad = \qquad  \max  & \quad \gamma \\
% \mbox{ subject to } &  \quad
%a^Tx - \gamma = \lmd f(x)+ \sum_{\nu\in\{0,1\}^m} \frac{g_1^{\nu_1}\cdots g_m^{\nu_m}}{p(x)}\sig_\nu(x),  \\
% & \qquad \qquad \lmd \geq 0,  \mbox{ every } \sig_\nu(x) \mbox{ is sos. }
%\end{align*}
%The dual of the above convex optimization problem is
%\begin{align*}
%  \min  \quad a^Tx  \qquad
%\mbox{ subject to }  \, \, x = \rho(y), \, y \in \mc{R}_{pre}^{\mc{D}}(f).
%\end{align*}
A similar contradiction argument like in (i) can be applied to prove the claim.
\end{proof}

\begin{example}
%(i)The convex set $\{x\in\re_+^2:\, x_1x_2^d \geq 1\}$ can be defined as
%$S_{\re_+^2}(f)$ with rational function
%$f(x) = x_2- \frac{1}{x_1^d}$.
%It can be verified that $f(x)$ is q-module concave over $\re_+^2$.
%So $S_{\re_+^2}(f)=\mc{R}_{qmod}^{\re_+^2}(f)$.\\
%(i) The convex set $\{(x_1,x_2)\in\re_+\times \re:\, x_1^2+x_1x_2 -x_2^2 \geq 1\}$ can be defined as
%$S_{\re_+\times \re}(f)$ with rational function
%$f(x) = x_1+x_2- \frac{x_2^2}{x_1}$.
%It can be verified that $f(x)$ is q-module concave over $\re_+^2$.
%So $S_{\re_+\times \re}(f)=\mc{R}_{qmod}^{\re_+\times \re}(f)$,
%which can be represented as
%\[
%\left\{(x_1,x_2):\, \exists \, z_0, z_1, z_2, s.t. \quad
%x_1 + x_2 - z_2 \geq 0, \,
%\bbm z_0 & 1 & z_1 \\ 1 & x_1 & x_2 \\ z_1 & x_2 & z_2 \ebm \succeq 0
%\right\}.
%\]
\begin{figure}
\centering
\includegraphics[width=.66\textwidth]{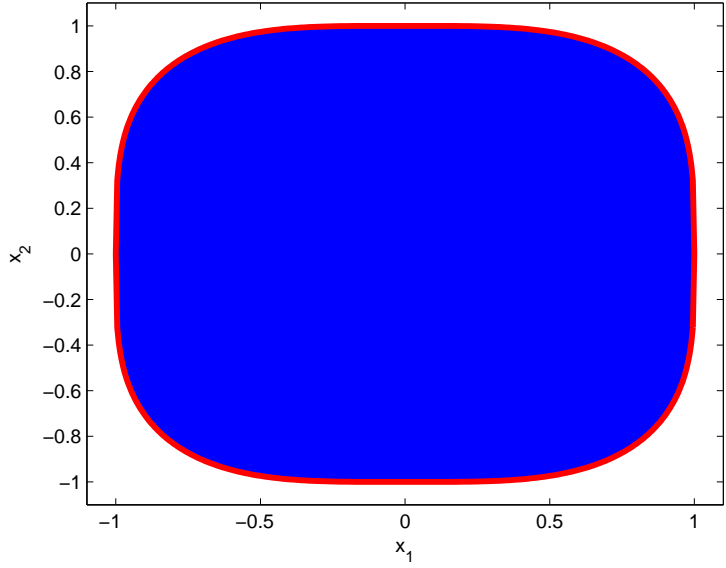}
\caption{ The convex set defined by
 $ x_1^2+x_2^2 \geq x_1^4+x_1^2x_2^2+x_2^4$.}
  \label{fig:RATx4}
\end{figure}
The convex set $\{x\in\re^2:\, x_1^2+x_2^2 \geq  x_1^4+x_1^2x_2^2+x_2^4 \}$ can be defined as
$S_{\re^2}(f)$ with rational function
$f(x) = 1- \frac{x_1^4+x_1^2x_2^2+x_2^4}{x_1^2+x_2^2}$.
The set $S_{\re^2}(f)$ is the shaded area bounded by a thick curve in Figure~\ref{fig:RATx4}.
We have already seen that $f(x)$ is first order sos concave.
So $S_{\re^2}(f)=\mc{R}_{qmod}^{\re^2}(f)$.
A polynomial division shows that
$\frac{1}{x_1^2+x_2^2}\mathfrak{m}_2(x)\mathfrak{m}_2(x)^T$ equals
{\small
\[
\bbm
0  &  0  &  0  &  1   &  0  &   0  \\
0  &  1  &  0  &  x_1 &  x_2 &  0 \\
0  &  0  &  0  &  x_2 &  0   &  0 \\
1  &  x_1 & x_2 & x_1^2-x_2^2  & x_1x_2 &  x_2^2  \\
0  &  x_2 & 0  &  x_1x_2 &  x_2^2  & 0  \\
0  &  0   &  0 &  x_2^2 &  0  &  0
\ebm  +  \frac{1}{p(x)}
\bbm
1  &  x_1  & x_2  & -x_2^2  & x_1x_2 &  x_2^2 \\
x_1  &  -x_2^2  &  x_1x_2 &  -x_1x_2^2  &  -x_2^3  & x_1x_2^2 \\
x_2  & x_1x_2 & x_2^2  & -x_2^3  &  x_1x_2^2  & x_2^3 \\
-x_2^2  & -x_1x_2^2 & -x_2^3  & x_2^4 & -x_1x_2^3  &  -x_2^4 \\
x_1x_2 & -x_2^3 & x_1x_2^2  & -x_1x_2^3  & -x_2^4  &  x_1x_2^3 \\
x_2^2 &  x_1x_2^2  & x_2^3  &  -x_2^4 & x_1x_2^3  & x_2^4
\ebm.
\]
}
So we can see that $S_{\re^2}(f)=\mc{R}_{qmod}^{\re^2}(f)$
can be represented as
\[
\left\{
\baray{l}
(x_1,x_2)\in \re\times \re: \qquad  \exists \quad y_{ij}, z_{ij}, \quad s.t. \qquad 1 \geq y_{20} + z_{04},  \\
\bbm
0  &  0  &  0  &  1   &  0  &   0  \\
0  &  1  &  0  &  x_1 &  x_2 &  0 \\
0  &  0  &  0  &  x_2 &  0   &  0 \\
1  &  x_1 & x_2 & y_{20}-y_{02}  & y_{11} &  y_{02}  \\
0  &  x_2 & 0  &  y_{11} &  y_{02}  & 0  \\
0  &  0   &  0 &  y_{02} &  0  &  0
\ebm  +
\bbm
z_{00}  &  z_{10}  & z_{01}  & -z_{02}  & z_{11} &  z_{02} \\
z_{10}  &  -z_{02}  &  z_{11} &  -z_{12}  &  -z_{03}  & z_{12} \\
z_{01}  & z_{11} & z_{02}  & -z_{03}  &  z_{12}  & z_{03} \\
-z_{02}  & -z_{12} & -z_{03}  & z_{04} & -z_{13}  &  -z_{04} \\
z_{11} & -z_{03} & z_{12}  & -z_{13}  & -z_{04}  &  z_{13} \\
z_{02} &  z_{12}  & z_{03}  &  -z_{04} & z_{13}  & z_{04}
\ebm \succeq 0
\earay \right\}.
\]
The plot of the projection of the above
coincides with the shaded area in Figure~\ref{fig:RATx4}.
\end{example}

%%Convexity of univariate rational functions

\iffalse

\subsubsection{Convexity of univariate rational functions}

\begin{example}
The univariate rational function $f(x)=\frac{1}{x^d}$ ($d>2$ is odd) is not globally convex,
but is convex over $\mc{D} = R_+$.
We are interested to know whether it is also q-module convex with respect to $\re_+$.
When $d=3$, the q-module convexity is obtained due to the identity
\[
\frac{1}{x^3}-\frac{1}{u^3} + \frac{3}{u^4} (x-u)  =  \frac{(x-u)^2}{x^3u^4}(2x^2+ (x-u)^2).
\]
When $d=5$, the q-module convexity can be observed by noting that
\begin{align*}
\frac{1}{x^5}-\frac{1}{u^5} + \frac{5}{u^6} (x-u)
& =  \frac{(x-u))^2}{x^5u^6}(u^4+2u^3x+3u^2x^2+4ux^3+5x^4)   \\
& =  \frac{(x-u)^2}{x^5u^6}\left( \half (x+u)^4 + (ux+x^2)^2 + \half (u^2-x^2)^2 + 3 x^4\right).
\end{align*}
Actually, we can show that $x^d$ is also always q-module convex with respect to $\re_+$
for any positive integer $d$.
\end{example}

\begin{pro}
Let $f(x)=\frac{q(x)}{p(x)}$ be a univariate rational functionc. Then we have
\bit
\item [(i)] If $f(x)$ is convex over $\re_+$,
then $f(x)$ must also be q-module convex over $\re_+$.

\item [(ii)] If $f(x)$ is convex over an interval $[a,b]$,
then $f(x)$ must also be q-module convex over $[a,b]$.
\eit
\end{pro}
\begin{proof}
???
\end{proof}
\fi

\subsection{Some special cases}

\subsubsection{Epigraph of rational functions}

The rational function $f(x)$ is convex over the convex domain $\mc{D}$
if and only if its epigraph
\[
\mathbf{epi}(f) := \{(x,t)\in \mc{D} \times \re: \,
f(x) \leq t \}
\]
is convex.
The LMI $\mc{R}_{qmod}^{\mc{D}}(t-f)$ and $\mc{R}_{pre}^{\mc{D}}(t-f)$
can be constructed by thinking of $t-f(x)$
as a polynomial in $x$ with coefficients in $t$.
So, if $f(x)$ is q-module (resp. preordering) convex over $\mc{D}$,
$\mc{R}_{qmod}^{\mc{D}}(t-f)$ (resp. $\mc{R}_{pre}^{\mc{D}}(t-f)$)
is an SDP representation for $\mathbf{epi}(f)$.

Now we consider the special case that
$f(x) = \frac{q(x)}{p(x)}$ is a univariate rational
function convex over an interval $I=[a,b]$
and $p(x)$ is positive over $(a,b)$.
Note that $I=\{x\in\re: x-a\geq 0, b-x\geq 0\}$
($x-a\geq 0$ is not required if $a=-\infty$, and similarly for $b-x\geq 0$).

\begin{theorem} \label{thm:epiURat}
Let $f(x) = \frac{q(x)}{p(x)}$ be a univariate rational
function and $p(x)$ is a polynomial nonnegative over an interval $I$.
If $f(x)$ is convex over $I$,
then its epigraph $\mathbf{epi}(f)=\mc{R}_{qmod}^{I}(f)$.
\end{theorem}

\begin{proof}
First assume $I=[a,b]$ is finite.
Since $\mc{R}_{qmod}^{[a,b]}(f)$ is linear in $t$,
it suffices to show for any fixed $t$
\[
\{x\in [a,b]:\, \exists y, z, s.t. , \,  x=y_1, (y,z,t)\in   \mc{R}_{qmod}^{[a,b]}(f) \}
= \{x\in [a,b]: f(x)\leq t \}.
\]
In the above the right hand side is contained in the left hand side.
Now we prove the converse.
For a contradiction, suppose there exists a tuple
$(\hat y, \hat z,t)\in   \mc{R}_{qmod}^{[a,b]}(f)$ such that
and $\hat x = \hat y_1$ does not belong to the convex set $\{x\in [a,b]: f(x)\leq t \}$.
Then there exists an affine function $c_0 + c_1 x$ such that
\[
c_0+c_1\hat x <0, \quad c_0+c_1 x \geq 0, \quad \forall \, x\in [a,b]: f(x) \leq t.
\]
Since $f(x)$ is convex over [a,b], by optimality condition,
there exists $\lmd \geq 0$ such that
\[
c_0+c_1 x -\lmd (t-f(x)) = \frac{(c_0+c_1 x -\lmd t) p(x) + \lmd q(x) }{p(x)} \geq 0,
\quad \forall \, x \in [a,b].
\]
Then we can see that $(c_0+c_1 x -\lmd t) p(x) + \lmd q(x)$ is a univariate polynomial
nonnegative on $[a,b]$.
So there exist sos polynomials $s_0(x), s_1(x), s_2(x)$
of degrees at most $2d,2d-2,2d-2$ respectively
(assume $\deg(p)+1, \deg(q) \leq 2d$ and see Powers and Reznick \cite{PR00})
such that
\[
(c_0+c_1 x -\lmd t) p(x) + \lmd q(x)  = s_0(x)+ (x-a)s_1(x)+(b-a)s_2(x)
\]
and hence then
\[
c_0+c_1 x = \lmd (t-f(x)) + \frac{ s_0(x)+ (x-a)s_1(x)+(b-a)s_2(x) }{p(x)}.
\]
Write $s_0(x), s_1(x), s_2(x)$ as
\[
s_0(x) = \mathfrak{m}_d(x)^T W_0 \mathfrak{m}_d(x), \quad
s_1(x) = \mathfrak{m}_{d-1}(x)^T W_1 \mathfrak{m}_{d-1}(x), \quad
s_2(x) = \mathfrak{m}_{d-1}(x)^T W_2 \mathfrak{m}_{d-1}(x)
\]
for some symmetric $W_0,W_1,W_2 \succeq 0$.
Then we have (noting $g_0=1, g_1 = x-a, g_2 = b-x$ and  \reff{R:mprod})
\begin{align*}
a^Tx - b  &=\lmd (t-f(x))+ \sum_{i=0}^2 
\left(\frac{g_i(x)}{p(x)} \mathfrak{m}_{d-d_i}(x)\mathfrak{m}_{d-d_i}(x)^T\right) \bullet W_i 
\quad (d_0=0, d_1=d_2 =1)\\
&= \lmd (t-f(x)) + 
\sum_{i=0}^2  \left(
\underset{\af \in \N: |\af| + |LE(p)| \leq 2d }{\sum}  Q_\af^{(i)}  x^\af
+ \underset{\bt \in \N: \bt < LE(p)}{\sum} P_\af^{(i)}  \frac{x^\bt}{p(x)} \right) \bullet W_i.
\end{align*}
In the above identity, 
if we replace each $x^\af$ by $\hat y_\af$
and $\frac{x^\bt}{p(x)}$ by $\hat z_\bt$,
then get the contradiction
\[
a^T \hat x - b = \lmd (t- (f^{(1)})^T\hat y - (f^{(2)})^T \hat z)   + 
\sum_{i=0}^2 Q_i(\hat y, \hat z) \bullet W_i \geq 0.
\]
Therefore we get
$\mathbf{epi}(f)= \mc{R}_{qmod}^{[a,b]}(f)$.

When $I$ is an infinite interval
of the form $[a,\infty)=\{x: x-a\geq 0\}$,
$(-\infty,b]=\{x: b-x\geq 0\}$ or
$(-\infty,\infty)=\{x: 1\geq 0\}$,
a similar argument can be applied to prove
$\mathbf{epi}(f)= \mc{R}_{qmod}^{I}(f)$.
\end{proof}

\subsubsection{Convex sets defined by structured rational functions}
\label{ssec:srtRAT}

For the convenience of discussion,
we define some basic convex sets ($r\leq s$)
\be  \label{def:Drt}
K_{r,s} = \{ (w,v)\in \re_+^r \times \re_+:\, w_1 \cdots w_r \geq v^s \}
\ee
which are all SDP representable (see \S3.3 in \cite{BTN}).

First, consider epigraphs of rational functions of the form
\be  \label{def:f}
f(x) =  \frac{\sum_{i=1}^P (q_i(x))^{s_i} }{p_1(x)\cdots p_r(x) }
\ee
where $q_i(x),p_j(x)$ are polynomials nonnegative over
$\mc{D}$ and the integer $s_i\geq r+1$. Then
\[
\mathbf{epi}(f) = \left\{ (x,t)\in \mc{D} \times \re:\,
\sum_{i=1}^P (q_i(x))^{s_i} \leq p_1(x)\cdots p_r(x) t
\right\}.
\]
The epigraph $\mathbf{epi}(f)$ can be equivalently presented as
\[
\mathbf{epi}(f) = \left\{ (x,t)\in \mc{D} \times \re:\,
\exists u_i, w_i,\,
\sum_{i=1}^P \frac{u_i^{s_i}}{w_1 \cdots w_r} \leq  t ,
q_i(x)\leq u_i, w_i \leq p_i(x), i=1,\ldots, P
\right\}.
\]
Note that
\[
q_i(x)\leq u_i  \Longleftrightarrow (x,u_i) \in \mathbf{epi}(q_i),
\qquad
w_i \leq p_i(x)  \Longleftrightarrow (x,-w_i) \in \mathbf{epi}(-p_i),
\]
\[
\sum_{i=1}^P \frac{u_i^{s_i}}{w_1 \cdots w_r} \leq  t  \Longleftrightarrow
\exists v_1, \cdots, v_P \geq 0, \,
v_1 + \cdots + v_P \leq t, (w_1,\cdots,w_r, v_i, u_i) \in K_{r+1,s_i}.
\]
Hence we obtain that
\begin{align*}  
\mathbf{epi}(f) = \Big\{ (x,t)\in \mc{D} \times \re:\,
\exists u_i, w_i, v_i, \,
\sum_{i=1}^P v_i \leq t, (w_1,\cdots,w_r, v_i, u_i) \in K_{r+1,s_i}, \\
(x,u_i) \in \mathbf{epi}(q_i),
(x,-w_i) \in \mathbf{epi}(-p_i), i=1,\ldots, P
\Big\}.
\end{align*}
When $q_i(x)$ are q-module (resp. preordering) convex over $\mc{D}$,
we know $\mathbf{epi}(q_i) = \mc{L}_{qmod}^{\mc{D}}(u_i-q_i)$
(resp. $\mathbf{epi}(q_i) = \mc{L}_{pre}^{\mc{D}}(u_i-q_i)$),
and similarly for $-p_i(x)$.
Thus, we get the theorem:

\begin{theorem} \label{thm:SRf}
Suppose $f(x)$ is given in the form \reff{def:f}
and all $q_i(x), p_j(x)$ there
are nonnegative over $\mc{D}$. 
\bdes
\item [(i)] If all $q_i(x), -p_j(x)$ are q-module convex over $\mc{D}$, then
\begin{align*}
\mathbf{epi}(f) = \Big\{ (x,t): \, \exists u_i, w_i, v_i,
x \in \bigcap_{i=1}^P \mc{L}_{qmod}^{\mc{D}}(u_i-q_i)
\bigcap_{j=1}^P \mc{L}_{qmod}^{\mc{D}}(p_j-w_j), \\
\sum_{i=1}^P v_i \leq t, (w_1,\cdots,w_r, v_i, u_i) \in K_{r+1,s_i}
\Big\}.
\end{align*}

\item [(ii)] If all $q_i(x), -p_j(x)$ are preordering convex over $\mc{D}$, then
\begin{align*}
\mathbf{epi}(f) = \Big\{ (x,t): \, \exists u_i, w_i, v_i,
x \in \bigcap_{i=1}^P \mc{L}_{pre}^{\mc{D}}(u_i-q_i)
\bigcap_{j=1}^P \mc{L}_{pre}^{\mc{D}}(p_j-w_j), \\
\sum_{i=1}^P v_i \leq t, (w_1,\cdots,w_r, v_i, u_i) \in K_{r+1,s_i}
\Big\}.
\end{align*}
\edes
\end{theorem}

\begin{example}
(i) Consider the epigraph
$\left\{(x_1,x_2,t)\in \re_+ \times \re \times \re: t \geq \frac{(1+x_2^2)^2}{x_1} \right\}.$
The rational function here is given in the form \reff{def:f}.
By Theorem~\ref{thm:SRf}, its epigraph can be represented as
\[
\left\{(x_1,x_2, t):\, \exists \, u, \,
\bbm t & u \\ u & x_1  \ebm \succeq 0, \,
\bbm u-1 & x_2 \\ x_2 & 1  \ebm \succeq 0
\right\}.
\]
(ii) Consider the epigraph
$
\left\{ (x,t)\in B(0,1) \times \re: \, t \geq \frac{(\sum_i x_i^2)^{n+1}}{(1-x_1^2)\cdots (1-x_n^2)} \right\}.
$
The rational function here is given in the form \reff{def:f}.
By Theorem~\ref{thm:SRf}, its epigraph can be represented as
\[
\left\{(x,t): \exists \,u, v_i, w_i, \, s.t., \, 
\sum_{i=1}^n v_i \leq t, \bbm 1-w_i & x_i \\ x_i & 1 \ebm \succeq 0, 
\bbm u & x^T  \\  x &  I_n  \ebm \succeq 0, 
(u,w_1,\ldots, w_n, u) \in K_{n+1,n+1}
\right\}.
\]
\end{example}

\medskip 
Second, consider epigraphs of rational functions of the form
\be \label{def:h}
h(x) = \sum_{k=1}^L  (f_k(x))^{b_k}
\ee
where $f_k(x)$ are rational functions given in form
\reff{def:f} and $b_k\geq 1$.
Then the epigraph $\mathbf{epi}(h)=\{(x,t)\in \mc{D}\times \re: h(x)\leq t\}$ can be presented as
\[
\mathbf{epi}(h) = \left\{ (x,t):\,
\exists \eta_k, \tau_k, \,\,\,\sum_{k=1}^L \tau_k \leq t,
(x,\eta_k) \in \mathbf{epi}(f_k), \,
(\tau_k, \eta_k) \in K_{1,b_k} , \,
k=1,\cdots,L
\right\}.
\]
Once the SDP representation for each $\mathbf{epi}(f_k)$ is available,
one SDP representation for $\mathbf{epi}(h)$ can be obtained 
consequently from the above.

\begin{example}
Consider the epigraph
$
\left\{(x_1,x_2,t)\in \re_+\times \re_+ \times \re:\,
t \geq \left(\frac{1+x_2^2}{x_1} \right)^2 \right\}.
$
From the above discussion, we know it can be represented as
\[
\left\{(x_1,x_2, t):\, \exists \, u, \,
\bbm t & u \\ u & 1 \ebm \succeq 0, \,
\bbm x_1 & 0 & 1 \\ 0 & x_1 & x_2 \\ 1& x_2 & u  \ebm \succeq 0
\right\}.
\]
\end{example}

\medskip

Third, consider the convex sets given in the form
\be \label{def:T}
T = \left\{ x\in \mc{D}:  a(x) \geq \sum_{k=1}^J h_k(x) \right\}
\ee
where $a(x)$ is a polynomial and every $h_k(x)$ is given in the form
\reff{def:h}.
Then
\[
T = \left\{ x:\, \exists t, \theta_1, \cdots, \theta_J, \,
t \geq \sum_{k=1}^J \theta_k, \,
(x,-t)\in \mathbf{epi}(-a(x)), (x,\theta_k)\in \mathbf{epi}(h_k(x)),
k=1,\ldots,J  \right\}.
\]
When $a(x)$ is q-module or preordering concave over $\mc{D}$,
$\mathbf{epi}(-a(x))$ is representable by
$\mc{L}_{qmod}^{\mc{D}}(a(x)-t)$ or $\mc{L}_{pre}^{\mc{D}}(a(x)-t)$.
Once the SDP representations for
$\mathbf{epi}(-a(x))$ and all $\mathbf{epi}(h_k(x))$ are available,
an SDP representation for $T$ can be obtained consequently.

\begin{example}
(i) Consider the convex set
$
T =\left\{(x_1,x_2)\in \re_+ \times \re:\, 1-3x_2^2 \geq \frac{(1+x_2^2)^2}{x_1} \right\}.
$
It is given in the form \reff{def:T} with $\mc{D} = \re_+ \times \re$.
From the above discussion, we have
\[ T=
\left\{(x_1,x_2):\, \exists \, u, w, \,
\bbm 1-u & x_2 \\ x_2 & \frac{1}{3}  \ebm \succeq 0, \,
\bbm w-1 & x_2 \\ x_2 & 1  \ebm \succeq 0, \,
(u,x_1,w) \in K_{2,2}
\right\}.
\]
Note that $(u,x_1,w) \in K_{2,2}$ has the representation
$\bbm u & w \\w & x_1  \ebm \succeq 0.$ \\
(ii) Consider the convex set
$
T =\left\{(x_1,x_2)\in \re_+\times \re:\, 1-x_2^2 \geq \frac{(1+x_2^2)^2}{x_1^2} \right\}.
$
It is given in the form \reff{def:T} with $\mc{D} = \re_+ \times \re$.
From the above discussion, we know it can be represented as
\[
\left\{(x_1,x_2):\, \exists \, u, \,
\bbm 1+u & x_2 \\ x_2 & 1-u \ebm \succeq 0, \,
\bbm x_1 & 0 & 1 \\ 0 & x_1 & x_2 \\ 1& x_2 & u  \ebm \succeq 0
\right\}.
\]
\end{example}

\section{Convex sets with singularities}
\label{sec:sig}

Let $S_\mc{D}(f) = \{x\in \mc{D}: f(x) \geq 0\}$
be a convex set defined by a polynomial or rational function $f(x)$.
Here $\mc{D} = \{x\in \re^n:\, g_1(x) \geq 0, \cdots, g_m(x)\geq 0\}$
is still a convex domain defined by polynomials.
Suppose the origin belongs to $S_\mc{D}(f)$
and is a singular point of the hypersurface $\mc{Z}(f)= \{x\in \cpx^n:\, f(x) = 0\}$, i.e.,
\[
f(0,0)=0,\, \nabla f(0,0) = 0.
\]
We are interested in finding 
SDP representability conditions for $S_\mc{D}(f)$.

As we have seen in Introduction,
one natural approach to getting an SDP representation for $S_\mc{D}(f)$
is to find a ``nicer'' defining function (possibly a concave rational function).
Let $p(x)$ be a polynomial or rational function positive in $int(\mc{D})$.
Then we can see $S_\mc{D}(f)$ is the closure of the set
\[
\left\{x\in int(\mc{D}): \frac{f(x)}{p(x)} \geq 0 \right\}.
\]
If $\frac{f(x)}{p(x)}$ has nice properties,
e.g., $\frac{f(x)}{p(x)}$ has special structures discussed in Section~\ref{sec:RatCvx},
or it is q-module or preordering concave over $\mc{D}$,
then an explicit SDP representation for $S_\mc{D}(f)$ can be obtained.
For instance, consider the convex set
\[
\{(x_1,x_2)\in \re_+^2:\, -x_1^3+3x_1x_2^2-(x_1^2+x_2^2)^2 \geq 0 \}.
\]
The origin is a singular point on its boundary.
If we choose $p(x) = x_2^2$, then it can be presented as
\[
\left\{(x_1, x_2)\in \re_+^2:\,
3x_1 \geq  2x_1^2 + x_2^2 + \frac{x_1^3}{x_2^2} + \frac{x_1^4}{x_2^2}
\right\}.
\]
Then this set can be represented as
\begin{align*}
\Big\{(x_1, x_2)\in \re_+^2:\, \exists \, u_1, u_2, u_3, u_4, \,\, s.t., \,
 3x_1 \geq 2u_1+u_2+u_3+u_4, \,
(u_1, x_1) \in K_{1,2}, \\
\qquad (u_2, x_2) \in K_{1,2},
(u_3, x_2,  x_2, x_1) \in K_{3,3},
(u_4, x_2,  x_2, x_1) \in K_{3,4}
\Big\}.
\end{align*}

However, there is no general procedure to find such a nice function $p(x)$.
In convex analysis, there is a technique called
{\it perspective transformation}
which might be very useful now.
Generally, we can assume
\[
S_\mc{D}(f) \subset \re_+ \times \re^{n-1}, \qquad
int(S_\mc{D}(f))  \ne \emptyset.
\]
%Note that the augmented set
%\[
%Aug(S_\mc{D}(f)) = \{(x_1,1,x_2,\cdots,x_n):  x\in S_\mc{D}(f) \}
%\]
%is also convex.
Define the perspective transformation $\mathfrak{p}$ as
\[
\mathfrak{p}(x_1,x_2,\cdots,x_n) = (1/x_1,x_2/x_1,\cdots,x_n/x_1).
\]
The image of $S_\mc{D}(f)$ under the perspective transformation $\mathfrak{p}$ is
\[
\Big\{\left(1/x_1,x_2/x_1,\cdots,x_n/x_1\right) :\,
x\in S_\mc{D}(f) \Big\} \subset \re_+ \times \re^{n-1},
\]
which is also convex (see \S2.3 in \cite{BV}). Define new coordinates
\[
\tilde x_1=\frac{1}{x_1}, \quad  \tilde x_2 = \frac{x_2}{x_1},  \quad  \cdots, \quad  \tilde x_n = \frac{x_n}{x_1}.
\]
Denote $\tilde{\tilde x} = ( \tilde x_2,  \cdots, \tilde x_n )$.
Suppose $f(x)$ has Laurent expansion around the origin
\[
f(x)=f_b(x) - f_{b+1}(x) - \cdots - f_d(x)
\]
where every $f_k(x)$ is a homogeneous part of degree $k$.
Let
\[
\tilde f( \tilde x):=\tilde f_0( \tilde{\tilde x} )
- \frac{\tilde f_{1}(\tilde{\tilde x})}{\tilde x_1} - \cdots -
\frac{\tilde f_{d-b}(\tilde{\tilde x} )}{\tilde x_1^{d-b}}
\]
where $\tilde f_i(\tilde x_2,\cdots,\tilde x_n ) := x_1^{b+i} f_{b+i}(\tilde x)$.
Define a new domain $\tilde{\mc{D}}$ as
\[
\tilde{\mc{D}} = \{\tilde x \in \re^n:\, \tilde g_1(\tilde x) \geq 0, \cdots, \tilde g_m( \tilde x) \geq 0 \}
\]
where $\tilde g_i(\tilde x)= \frac{g_i(x)}{x_1^{\deg(g_i)}}$.
Note that $\tilde{\mc{D}}$ is convex if and only if $\mc{D}$ is convex (see \S2.3 in \cite{BV}).
Therefore, under the perspective transformation $\mathfrak{p}$,
the set $S_\mc{D}(f)$ can be equivalently defined as
\[
S_{\tilde{\mc{D}}}(\tilde f)  = \{\tilde x\in \tilde{\mc{D}}: \tilde f( \tilde x )\geq 0 \}.
\]

\begin{pro} \label{pro:h0psd}
If $S_{\mc{D}}(f)$ is convex,
then $\tilde f_0(\tilde{\tilde x}) \geq 0$ for any
$\tilde x = (\tilde x_1, \tilde{\tilde x}) \in S_{\tilde{\mc{D}}}(\tilde f)$.
\end{pro}
\begin{proof}
Fix $x = (x_1,\cdots, x_n ) \in S_{\mc{D}}(f)$ and
$\tilde x = (\tilde x_1,\cdots,\tilde x_n ) \in S_{\tilde{\mc{D}}}(\tilde f)$
such that $\tilde x = \mathfrak{p}(x_1,x_2,\cdots,x_n)$.
By the convexity of $S_{\mc{D}}(f)$, the line segment
$\{tx: 0\leq t\leq 1\}$ belongs to $S_{\mc{D}}(f)$.
Thus its image
\[
\mathfrak{p}(tx_1,tx_2,\cdots,tx_n) = (\frac 1t \tilde x_1, \tilde x_2, \cdots,\tilde x_n )
\]
belongs to $S_{\tilde{\mc{D}}}(\tilde f)$.
Now let $t\to 0$. Then $ \tilde f (\frac 1t \tilde x_1, \tilde x_2, \cdots,\tilde x_n ) \geq 0 $
implies $\tilde f_0( \tilde{\tilde x} )\geq 0$.
\end{proof}

\subsection{ The case of structured $\tilde f_k (\tilde{\tilde x})$ }

In this subsection, we assume $\tilde f_k (\tilde{\tilde x})$ have special structures.
Then the methods in Subsection~\ref{ssec:srtRAT}
can be applied to construct SDP representations for
$S_{\tilde{\mc{D}}}(\tilde f)$.

\begin{theorem} \label{thm:strSig}
Suppose every $\tilde f_k (\tilde{\tilde x})\, (k=1,\cdots,d-b)$ is given in the form
\[
\tilde f_k (\tilde{\tilde x}) = \sum_{i=1}^{L_k} (q_{k,i}(\tilde{\tilde x}) )^{r_{k,i}} +
\sum_{j=1}^{R_k} \Big( \sum_{\ell = 1}^{Q_{k,j}} (p_{k,j,\ell}(\tilde{\tilde x}))^{s_{k,j,\ell} } \Big)^{k}
\]
for some polynomials $q_{k,i}(\tilde{\tilde x}), p_{k,j,\ell}(\tilde{\tilde x})$ 
which are nonnegative over $\tilde{\mc{D}}$
and integers $r_{k,i} \geq k+1$, $s_{k,j,\ell} \geq 2$
($p_{k,j}(\tilde{\tilde x})$ can be any affine polynomial when $s_{k,j,\ell}$ is even).
Then $S_{\tilde{\mc{D}}}(\tilde f)$
can be represented as
\begin{align*}
\left\{ \tilde x :\, \exists \, u, u_{k,i}, v_{k,j,\ell}, \,
 u \geq \sum_{k=1}^{d-b}
\left(\sum_{i=1}^{L_k} \frac{(u_{k,i})^{r_{k,i}}}{\tilde x_1^k} +
\sum_{j=1}^{R_k} 
\left( \sum_{\ell = 1}^{Q_{k,j}} \frac{(v_{k,j,\ell})^{s_{k,j,\ell}}}{\tilde x_1}
 \right)^{k} \right)\right.,\\
(\tilde{\tilde x}, -u) \in \mathbf{epi}( -\tilde f_0 ), \,
(\tilde{\tilde x}, u_{k,i}) \in \mathbf{epi}( q_{k,i} ), \,
(\tilde{\tilde x}, v_{k,j,\ell}) \in \mathbf{epi}( p_{k,j,\ell} )
\Big\}.
\end{align*}
Furthermore, if $-\tilde f_0(\tilde{\tilde x})$ and all
$q_{k,j}(\tilde{\tilde x}), p_{k,j,\ell}(\tilde{\tilde x})$ are q-module or
preordering convex over $\tilde{\mc{D}}$,
then $S_{\tilde{\mc{D}}}(\tilde f)$ is SDP representable.
\end{theorem}
\begin{proof}
The first conclusion is obvious by introducing new variables
$u, u_{k,i}, v_{k,j,\ell}$.
Note that 
\[
u \geq \sum_{k=1}^{d-b}
\Big(\sum_{i=1}^{L_k} \frac{(u_{k,i})^{r_{k,i}}}{\tilde x_1^k} +
\sum_{j=1}^{R_k} \Big( \sum_{\ell = 1}^{Q_{k,j}} \frac{(v_{k,j,\ell})^{s_{k,j,\ell}}}{\tilde x_1} \Big)^{k} \Big)
\]
is equivalent to
\begin{align*}
u \geq \sum_{k=1}^{d-b}
\Big(\sum_{i=1}^{L_k} \eta_{k,i}   +
\sum_{j=1}^{R_k} \xi_k \Big),
(\tilde x_1,\cdots,\tilde x_1, \,\,  \eta_{k,i}, u_{k,i}) &\in K_{k+1,r_{k,i}}, \\
(\xi_k, \zeta_k) \in K_{1,k}, \,\, (\zeta_k,\tilde x_1, v_{k,j,\ell})& \in K_{2,s_{k,j,\ell}}.
\end{align*}
When $-\tilde f_0(\tilde{\tilde x})$
is q-module (resp. preordering) convex over $\tilde{\mc{D}}$,
$(\tilde{\tilde x}, -u) \in \mathbf{epi}( -\tilde f_0 )$
is representable by $\mc{L}_{qmod}^{\tilde{\mc{D}}}(u-\tilde f_0)$
(resp. $\mc{L}_{pre}^{\tilde{\mc{D}}}(u-\tilde f_0)$).
Similar results hold for $q_{k,j}(\tilde{\tilde x}), p_{k,j,\ell}(\tilde{\tilde x})$.
Once the SDP representations for epigraphs of $-\tilde f_0(\tilde{\tilde x})$,
$q_{k,j}(\tilde{\tilde x})$ and $p_{k,j,\ell}(\tilde{\tilde x})$
are all available, we can get an SDP representation for
$S_{\tilde{\mc{D}}}(\tilde f)$ consequently.
\end{proof}

%Obviously, $S_{\tilde{\mc{D}}}(\tilde f)$
%can be represented as
%\[
%\left\{ \tilde x \in \tilde{\mc{D}} : \, \exists \, t, t_1, \ldots, t_{d-b}, \,
%t\geq t_1 + \cdots + t_{d-b}, \,\,
%\tilde f_0 (\tilde{\tilde x}) \geq t, \,\,
%t_k \geq \frac{\tilde f_k (\tilde{\tilde x}) }{\tilde x_1^k},
%k=1,\ldots,d-b
%\right\}.
%\]
%When $\tilde f_0 (\tilde{\tilde x})$ is q-module concave (resp. preodering concave)
%over $\tilde{\mc{D}}$, $\tilde f_0 (\tilde{\tilde x}) \geq t$ can be represented
%by $\mc{L}_{qmod}^{\tilde{\mc{D}}}(\tilde f_0 - t)$
%(resp. $\mc{L}_{pre}^{\tilde{\mc{D}}}(\tilde f_0 -t)$).
%Once all the SDP representations of
%$\tilde f_0 (\tilde{\tilde x}) \geq t$ and $t_k \geq \frac{\tilde f_k (\tilde{\tilde x}) }{\tilde x_1^k}$
%are obtained,
%we can get one for $S_{\tilde{\mc{D}}}(\tilde f)$ immediately.
%

Now we show some examples on how to apply Theorem~\ref{thm:strSig}.

%%%%%%%%%%%%%%%
\begin{figure}[htb]
\centering
\btab{cc}
\includegraphics[width=.45\textwidth]{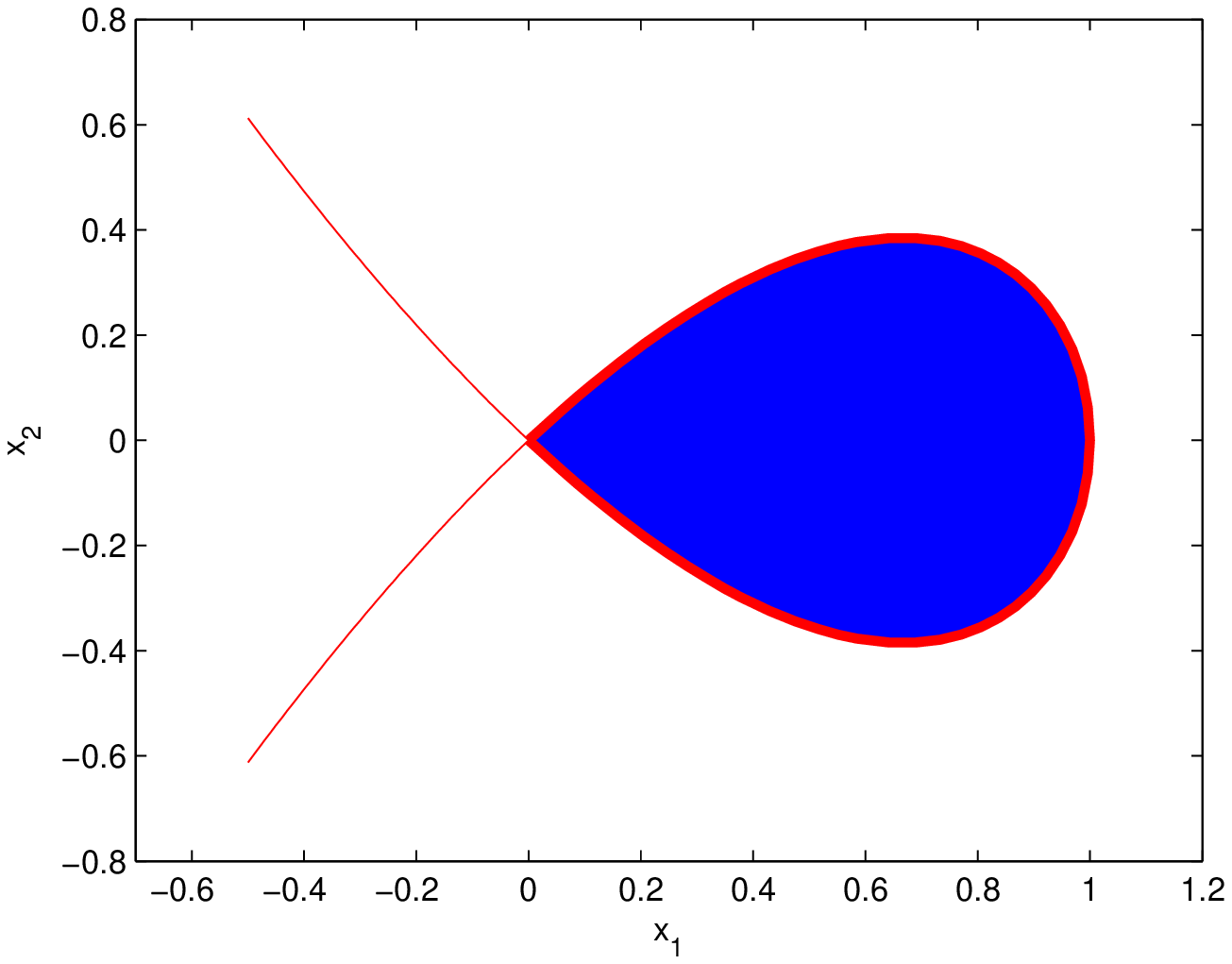} &
\includegraphics[width=.45\textwidth]{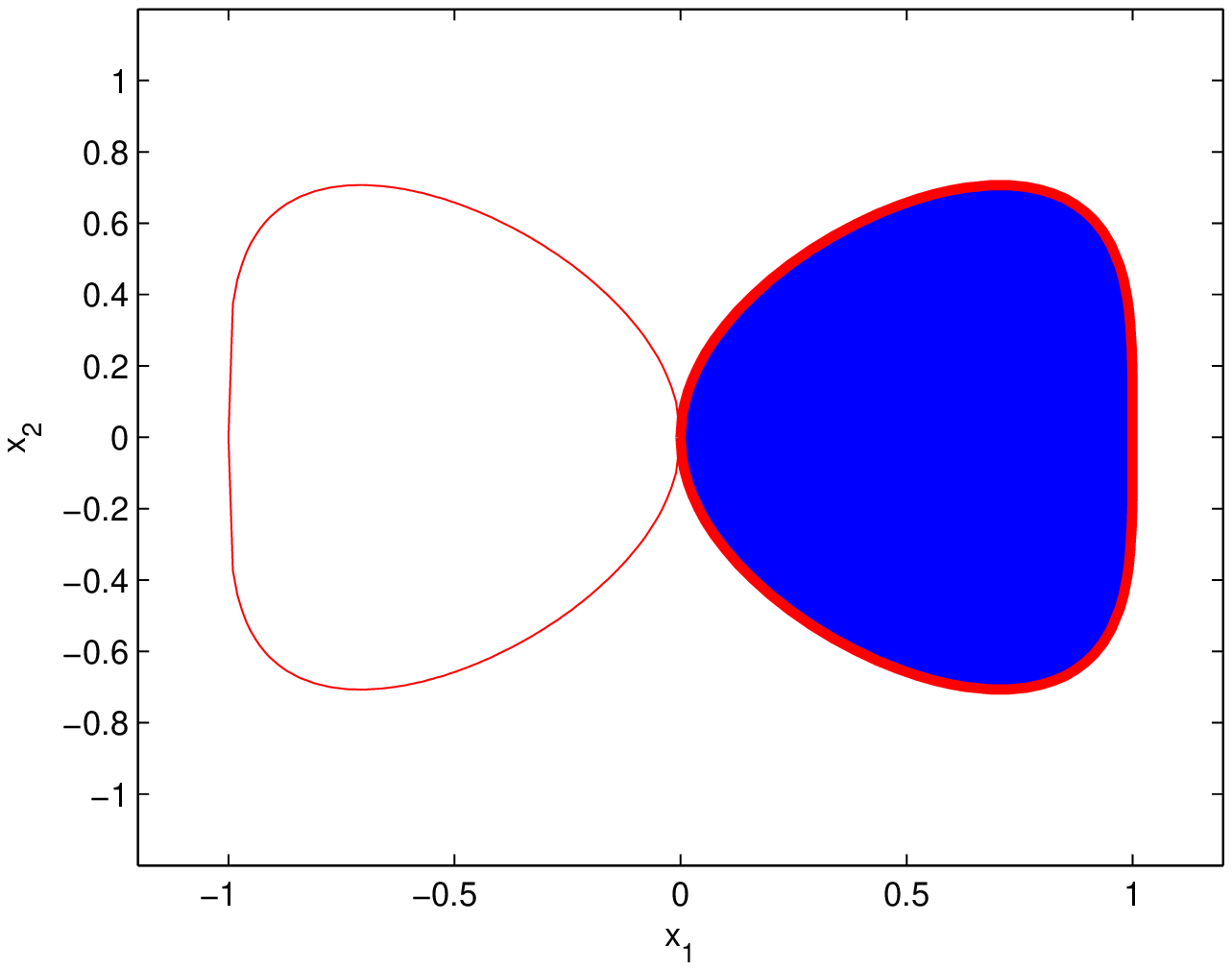} \\
(a): $x_1^2-x_1^3-x_2^2 \geq 0, x_1 \geq 0$ &
(b): $x_1^2-x_1^4-x_2^4 \geq 0, x_1 \geq 0$ \\
\includegraphics[width=.45\textwidth]{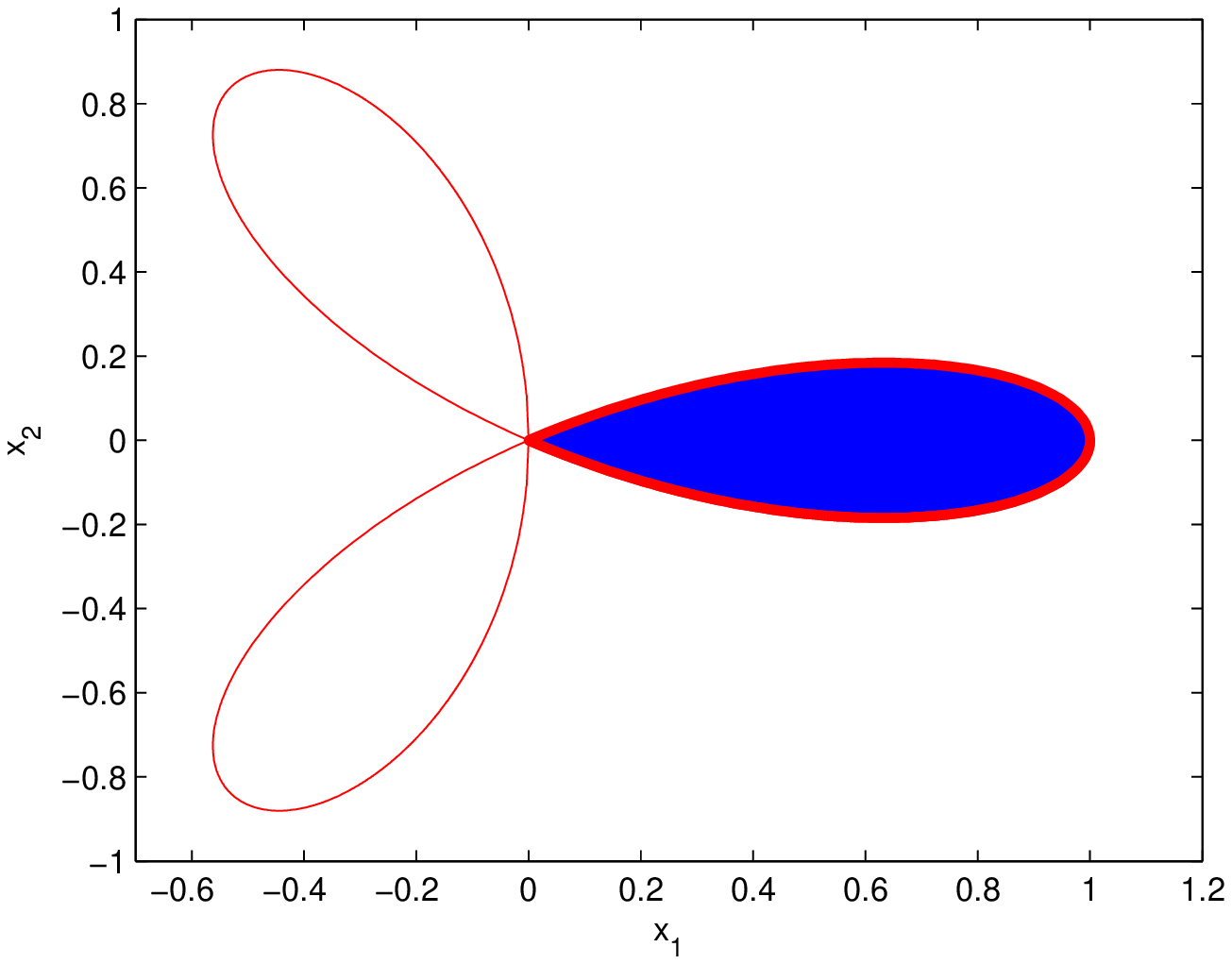} &
\includegraphics[width=.45\textwidth]{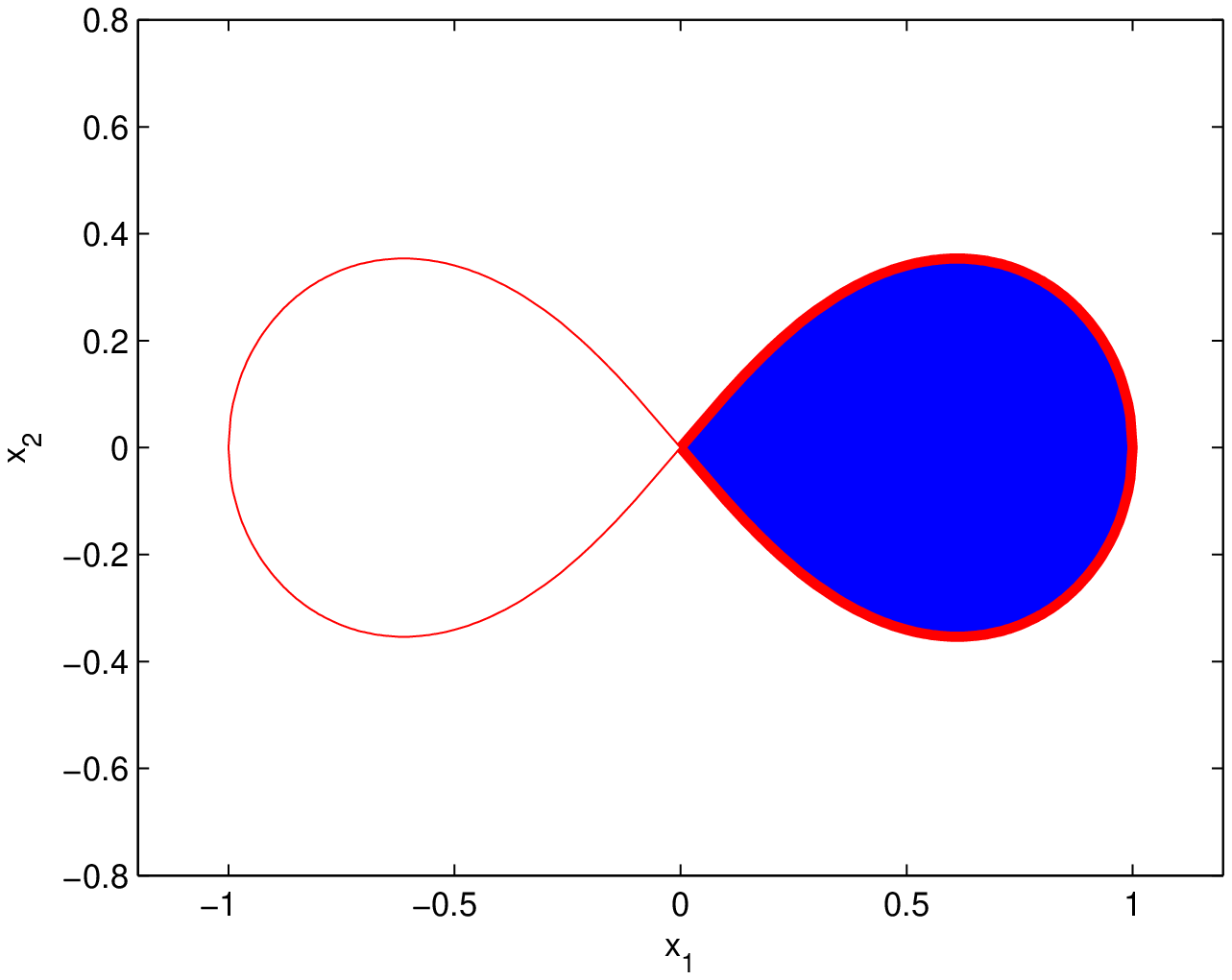} \\
(c): $x_1^3-3x_1x_2^2-(x_1^2+x_2^2)^2 \geq 0, x_1 \geq 0$ &
(d): $x_1^2-x_2^2-(x_1^2+x_2^2)^2 \geq 0, x_1 \geq 0$
\etab
\caption{Four singular convex sets discussed in Example~\ref{exmp:4sig}}.
\label{4SigCurv}
\end{figure}
%%%%%%%%%%%%%%%%%%%%%%%%%%%%%%%%%%

\begin{example}  \label{exmp:4sig}
(a) Consider convex set
$
S =\{(x_1,x_2)\in \re_+\times \re:\, x_1^2-x_1^3-x_2^2\geq 0 \}.
$
Its boundary is a cubic curve and
the origin is a singular node.
% \begin{figure}
%\centering
%\includegraphics[width=.66\textwidth]{pic_node.eps}
%\caption{ The singular convex set defined by
% $x_1^3-3x_1x_2^2-(x_1^2+x_2^2)^2 \geq 0, x_1 \geq 0$ .}
% \label{fig:node}
%\end{figure}
This convex set is the shaded area
bounded by a thick curve in Figure~\ref{4SigCurv}(a).
The thin curves are other branches of this cubic curve.
After the perspective transformation, we get
\[
\tilde S =\left\{(\tilde x_1, \tilde x_2)\in \re_+ \times \re:\,
1-\tilde x_2^2-\frac{1}{\tilde x_1} \geq 0\right\}
\]
which can be represented as
\[
\left\{(\tilde x_1, \tilde x_2):\, \exists \, u, \,\, s.t., \,
\bbm 1- u & \tilde x_2 \\ \tilde x_2 & 1 \ebm \succeq 0, \,
\bbm u & 1 \\ 1 & \tilde x_1 \ebm \succeq 0
\right\}.
\]
After the inverse perspective transformation, we get
\[
S =
\left\{( x_1,   x_2):\, \exists \, u, \,\, s.t., \,
\bbm x_1- u &  x_2 \\   x_2 & x_11 \ebm \succeq 0, \,
\bbm u & x_1 \\ x_1 & x_1 \ebm \succeq 0
\right\}.
\]
The plot of the projection of the above
coincides with the shaded area in Figure~\ref{4SigCurv}(a). \\
(b) Consider convex set
$
S =\{(x_1,x_2)\in \re_+\times \re:\, x_1^2-x_2^4-x_2^4 \geq 0 \}.
$
The origin is a singular tacnode on the boundary.
% \begin{figure}
%\centering
%\includegraphics[width=.66\textwidth]{pic_tacnode.eps}
%\caption{ The singular convex set defined by
% $x_1^2-x_2^4-x_2^4 \geq 0, x_1 \geq 0$ .}
% \label{fig:tacnode}
%\end{figure}
This convex set is the shaded area
bounded by a thick curve in Figure~\ref{4SigCurv}(b).
The thin curve is the other branch of the singular curve
$x_1^2-x_2^4-x_2^4 = 0$.
After the perspective transformation, we get
\[
\tilde S =\left\{(\tilde x_1, \tilde x_2)\in \re_+ \times \re:\,
1 \geq \frac{1+\tilde x_2^4}{\tilde x_1^2}\right\}
\]
which can be represented as
\[
\left\{(\tilde x_1, \tilde x_2):\, \exists \, u_1,u_2, \,\, s.t., \,
\bbm 1- u_1 & u_2 \\ u_2 & 1+u_1 \ebm \succeq 0, \,
\bbm u_1 & 1 \\ 1 & \tilde x_1 \ebm \succeq 0, \,
\bbm u_2 & \tilde x_2 \\ \tilde x_2  &\tilde x_1\ebm \succeq 0
\right\}.
\]
After the inverse perspective transformation, we get
\[
S = \left\{( x_1, x_2):\, \exists \, u_1,u_2, \,\, s.t., \,
\bbm x_1- u_1 & u_2 \\ u_2 & x_1+u_1 \ebm \succeq 0, \,
\bbm u_1 & x_1 \\ x_1 & 1 \ebm \succeq 0, \,
\bbm u_2 & x_2 \\ x_2  & 1\ebm \succeq 0
\right\}.
\]
The plot of the projection of the above
coincides with the shaded area in Figure~\ref{4SigCurv}(b). \\
(c) Consider convex set
$
S =\{(x_1,x_2)\in \re_+\times \re:\, x_1^3-3x_1x_2^2-(x_1^2+x_2^2)^2 \geq 0 \}.
$
The origin is a singular point on the boundary.
This convex set is the shaded area
bounded by a thick curve in Figure~\ref{4SigCurv}(c).
The thin curves are other branches of the singular curve
$x_1^3-3x_1x_2^2-(x_1^2+x_2^2)^2 = 0$.
After the perspective transformation, we get
\[
\tilde S =\left\{(\tilde x_1, \tilde x_2)\in \re_+ \times \re:\,
1-3\tilde x_2^2-\frac{(1+\tilde x_2^2)^2}{\tilde x_1} \geq 0\right\}
\]
which can be represented as
\[
\left\{(\tilde x_1, \tilde x_2):\, \exists \, u, w, \,\, s.t., \,
\bbm 1- u & \tilde x_2 \\ \tilde x_2 & \frac{1}{3} \ebm \succeq 0, \,
\bbm u & w \\ w & \tilde x_1 \ebm \succeq 0, \,
\bbm w-1 & \tilde x_2 \\ \tilde x_2 & 1 \ebm \succeq 0
\right\}.
\]
After the inverse perspective transformation, we get
\[
S = \left\{(x_1,x_2):\, \exists\, u,w,\, s.t.,\,
\bbm x_1- u & x_2 \\  x_2 & \frac{1}{3} x_1 \ebm \succeq 0, \,
\bbm u & w \\ w & 1 \ebm \succeq 0, \,
\bbm w-x_1 &   x_2 \\   x_2 & x_1 \ebm \succeq 0
\right\}.
\]
The plot of the projection of the above
coincides with the shaded area in Figure~\ref{4SigCurv}(c). \\
(d) Consider convex set
$
S =\{(x_1,x_2) \in \re_+\times \re:\, x_1^2-x_2^2-(x_1^2+x_2^2)^2 \geq 0 \}.
$
The origin is a singular point on the boundary
which is one branch of the lemniscate curve
$x_1^2-x_2^2-(x_1^2+x_2^2)^2 = 0$.
%\begin{figure}
%\centering
%\includegraphics[width=.66\textwidth]{pic_lemni.eps}
%\caption{ The singular convex set defined by
% $x_1^2-x_2^2-(x_1^2+x_2^2)^2 \geq 0, x_1 \geq 0$ .}
%\label{fig:lemni}
%\end{figure}
This convex set is the shaded area bounded
a thick lemniscate curve in Figure~\ref{4SigCurv}(d).
The thin curve is the other branch of the lemniscate curve.
After the perspective transformation, we get
\[
\tilde S =\left\{(\tilde x_1,\tilde x_2)\in \re_+\times \re:\,
(1-\tilde x_2^2) \geq \left(\frac{1+\tilde x_2^2}{\tilde x_1}\right)^2 \right\}
\]
which can be represented as
\[
\left\{(\tilde x_1, \tilde x_2):\, \exists \, u, w, \,\, s.t. \,
\bbm 1+ u & \tilde x_2 \\ \tilde x_2 & 1-u \ebm \succeq 0, \,
\bbm u & 1 & \tilde x_2 \\ 1 & \tilde x_1 & 0 \\ \tilde x_2 & 0 & \tilde x_1  \ebm \succeq 0
\right\}.
\]
After the inverse perspective transformation, we get
\[
S =
\left\{(x_1, x_2):\, \exists \, u, \,\, s.t. \,
\bbm x_1+ u & x_2 \\  x_2 & x_1 -u \ebm \succeq 0, \,
\bbm u & x_1 &  x_2 \\ x_1 &  1 & 0 \\   x_2 & 0 & 1  \ebm \succeq 0
\right\}.
\]
The plot of the projection of the above
coincides with the shaded area in Figure~\ref{4SigCurv}(d).
\end{example}

\subsection{The case of two consecutive homogeneous parts}

In this subsection, we consider the special case that
$f(x) = f_b(x) -f_{b+1}(x)$ having two consecutive homogeneous parts.
Then, after perspective transformation, we get
\[
S_{\tilde{\mc{D}} }(\tilde f ) = \left\{ \tilde x \in \tilde{\mc{D}}:
\tilde f_0(\tilde{\tilde x} ) -
\frac{\tilde f_1(\tilde{\tilde x})}{\tilde x_1 } \geq 0
\right\}.
\]
By Proposition~\ref{pro:h0psd}, for any $\tilde x \in \tilde{\mc{D}}$,
we have $f_0(\tilde{\tilde x} )\geq 0$.
Let $\mc{D}^\prm$ be the intersection of
$\{\tilde{\tilde x} : f_0( \tilde{\tilde x}  )\geq 0\}$
and the projection of $ \mc{D} \subset \re_+ \times \re^{n-1}$
into $\re^{n-1}$. Then we get
\[
S_{\tilde{\mc{D}} }(\tilde f ) =
\left\{ (\tilde{\tilde x},  \tilde x_1)  \in \mc{D}^\prm \times \re_+:
\tilde x_1 \geq h(\tilde{\tilde x})
\right\}, \quad
h(\tilde{\tilde x}):=
\frac{\tilde f_1(\tilde{\tilde x})}{\tilde f_0(\tilde{\tilde x})}.
\]
Then $S_{\tilde{\mc{D}} }(\tilde f )$ is
the epigraph $\mathbf{epi}(h)$ of the rational function $h(\tilde{\tilde x})$
over $\mc{D}^\prm$.
Note that $h(\tilde{\tilde x})$ is convex over the domain $\mc{D}^\prm$
if and only if $S_{\tilde{\mc{D}} }(\tilde f )$ is convex,
which then holds if and only if $S_{\mc{D} }( f )$ is convex.
So, if $h(\tilde{\tilde x})$ is q-module or preordering convex over $\mc{D}^\prm$,
then $\mc{R}_{qmod}^{\mc{D}^\prm}(\tilde x_1 - h)$ or $\mc{R}_{pre}^{\mc{D}^\prm}(\tilde x_1 -h)$
is an SDP representation for $\mathbf{epi}(h)$,
and then one can be obtained for $S_{\mc{D}}(h)$
after the inverse perspective transformation.

A very interesting case is $n=2$.
Then $\mc{D}^\prm$ must be an interval $I$ of the real line.

\begin{theorem} \label{thm:2D2Tsig}
Let $n=2$ and $\mc{D}^\prm = I$ be an interval as above.
If $S_\mc{D}(f)$ is convex, then
$S_{\tilde{\mc{D}} }(\tilde f )= \mc{R}_{qmod}^{I}(\tilde x_1 - h)$,
and hence $S_\mc{D}(f)= \mathfrak{p}^{-1}(\mc{R}_{qmod}^{I}(\tilde x_1 - h))$.
\end{theorem}
\begin{proof}
When $n=2$, $f(\tilde{\tilde x} )$ is a univariate rational function.
When $S_{\mc{D} }(f)$ is convex,
$S_{\tilde{\mc{D}} }(\tilde f)$ is also convex.
Since $S_{\tilde{\mc{D}} }(\tilde f)$ is the epigraph of $f(\tilde{\tilde x} )$,
$f(\tilde{\tilde x} )$ must be a univariate rational function
convex over the interval $I$.
By Theorem~\ref{thm:epiURat}, its epigraph is representable by
$\mc{R}_{qmod}^{I}(\tilde x_1 - h)$.
Thus $S_{\tilde{\mc{D}} }(\tilde f)= \mc{R}_{qmod}^{I}(\tilde x_1 - h)$.
After the inverse perspective transformation,
we can get an SDP representation for $S_{\mc{D} }(f)$.
\end{proof}

Now we see some examples on how to
find SDP representations for singular convex sets
by applying Theorem~\ref{thm:2D2Tsig}.

\begin{example}
(i)
We revisit the singular convex set (a) in Example~\ref{exmp:4sig}.
After the perspective transformation, we get
$
\tilde S =\left\{(\tilde x_1, \tilde x_2) \in  \re \times [-1, 1]:\,
\tilde x_1 \geq  \frac{1}{ 1- \tilde x_2^2 }\right\},
$
which can be represented as
\[
\left\{(\tilde x_1, \tilde x_2):
\baray{lr}
\exists \, y_2, y_3, y_4, z_0, z_1,  &
\tilde x_1 \geq z_0, \qquad
\bbm 1 - y_2 & \tilde x_2 - y_3 \\ \tilde x_2 - y_3 & y_2 - y_4 \ebm \succeq 0 , \\
\bbm 1  & \tilde x_2 & y_2 \\  \tilde x_2 & y_2 & y_3 \\  y_2 & y_3 &  y_4 \ebm \succeq 0 , &
\bbm z_0 & z_1  & -1+z_0 \\
z_1 & -1+z_0  & -\tilde x_2 + z_1  \\
-1+z_0  & -\tilde x_2 + z_1  & - y_2 -1 + z_0
\ebm \succeq 0
\earay
\right\}.
\]
Applying the inverse perspective transformation, we get an SDP representation for $S$
\[
\left\{(x_1, x_2):
\baray{lr}
\exists \, y_2, y_3, y_4, z_0, z_1, &
1 \geq z_0, \qquad
\bbm x_1 - y_2 &   x_2 - y_3 \\ x_2 - y_3 & y_2 - y_4 \ebm \succeq 0 , \\
\bbm x_1  &   x_2 & y_2 \\   x_2 & y_2 & y_3 \\  y_2 & y_3 &  y_4 \ebm \succeq 0 , &
\bbm z_0 & z_1  & -x_1+z_0 \\
z_1 & -x_1+z_0  & -x_2 + z_1  \\
-x_1+z_0  & -x_2 + z_1  & - y_2 -x_1 + z_0
\ebm \succeq 0
\earay
\right\}.
\]
Interestingly but not surprisingly,
the plot of the above coincides with the shaded area in Figure~\ref{4SigCurv}(a).  \\
(ii) Revisit the singular convex set (c) in Example~\ref{exmp:4sig}.
After the perspective transformation, we get
$
\tilde S =\left\{(\tilde x_1, \tilde x_2) \in \re \times \frac{1}{\sqrt{3}}[-1, 1]:\,
\tilde x_1 \geq -\frac{x_2^2}{3} -\frac{7}{9} + \frac{16}{ 9(1-3x_2^2)}\right\},
$
which equals
\[
\left\{(\tilde x_1, \tilde x_2):
\baray{lr}
\exists \, y_2, y_3, y_4, z_0, z_1,  &
\tilde x_1 \geq -\frac{1}{3} y_2 -\frac{7}{9} + \frac{16 }{9} z_0,
\bbm 1 - 3 y_2 & \tilde x_2 - 3 y_3 \\ \tilde x_2 - 3 y_3 & y_2 - 3 y_4 \ebm \succeq 0 , \\
\bbm 1  & \tilde x_2 & y_2 \\  \tilde x_2 & y_2 & y_3 \\  y_2 & y_3 &  y_4 \ebm \succeq 0 , &
\bbm z_0 & z_1  & \frac{1}{3}(z_0-1) \\
z_1 & \frac{1}{3}(z_0-1)  & \frac{1}{3}(z_1-\tilde x_2)  \\
\frac{1}{3}(z_0-1)  & \frac{1}{3}(z_1-\tilde x_2) & -\frac{1}{3} y_2 + \frac{1}{9}(z_0-1)
\ebm \succeq 0
\earay
\right\}.
\]
Applying the inverse perspective transformation, we get an SDP representation for $S$
\[
\left\{(x_1, x_2):
\baray{lr}
\exists \, y_2, y_3, y_4, z_0, z_1, &
1 \geq -\frac{1}{3} y_2 -\frac{7}{9} x_1 + \frac{16 }{9} z_0,
\bbm x_1 - 3 y_2 & x_2 - 3 y_3 \\ x_2 - 3 y_3 & y_2 - 3 y_4 \ebm \succeq 0,  \\
\bbm x_1 & x_2 & y_2 \\  x_2 & y_2 & y_3 \\  y_2 & y_3 &  y_4 \ebm \succeq 0, &
\bbm
z_0 & z_1  & \frac{1}{3}(z_0-x_1) \\
z_1 & \frac{1}{3}(z_0-x_1)  & \frac{1}{3}(z_1- x_2)  \\
\frac{1}{3}(z_0-x_1)  & \frac{1}{3}(z_1- x_2) & -\frac{1}{3} y_2 + \frac{1}{9}(z_0-x_1)
\ebm \succeq 0
\earay
\right\}.
\]
Also interestingly but not surprisingly,
the plot of the above also
coincides with the shaded area in Figure~\ref{4SigCurv}(c).
\end{example}

Now let us conclude this subsection with an example
such that Theorem~\ref{thm:2D2Tsig}
can be applied to get an SDP representation for $S_{\mc{D}}(f)$
while Theorem~\ref{thm:strSig} can not.

\begin{figure}[htb]
\centering
\includegraphics[width=.66\textwidth]{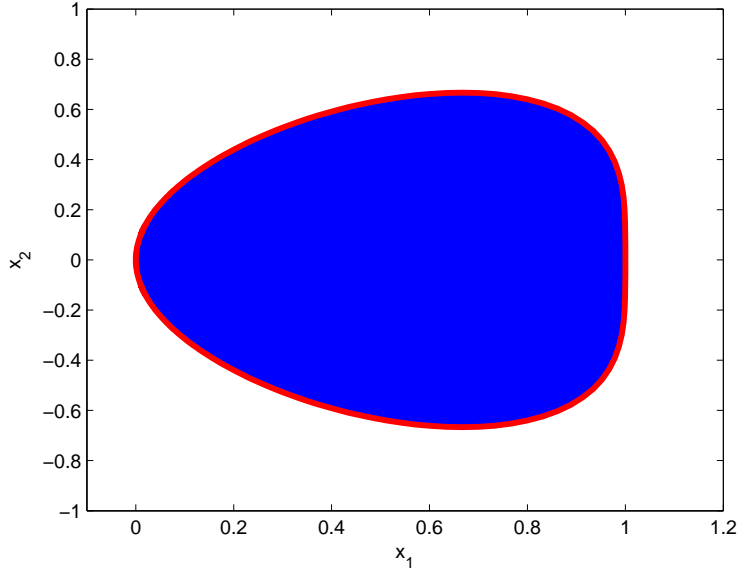}
\caption{ The singular convex set defined by
 $x_1(x_1^2+x_2^2)- x_1^4-x_1^2x_2^2-x_2^4 \geq 0 $ .}  \label{fig:bean}
\end{figure}

\begin{example}
Consider the convex set
$
S =\{(x_1,x_2) \in \re^2:\, x_1(x_1^2+x_2^2)- x_1^4-x_1^2x_2^2-x_2^4 \geq 0 \}.
$
The origin is a singular point on the boundary $\pt S$
which is a quartic bean curve.
The picture of this convex set is the shaded area
bounded by the thick bean curve in Figure~\ref{fig:bean}.
After the perspective transformation, we get
$
\tilde S =\left\{(\tilde x_1,\tilde x_2)\in \re^2:\,
\tilde x_1 \geq f(\tilde x_2):= \tilde x_2^2 + \frac{1}{1+\tilde x_2^2} \right\}.
$
$f(\tilde x_2)$ does not have structures required by Theorem~\ref{thm:strSig}.
Obviously $\mc{D}^\prm =  (-\infty,\infty)$.
We can check that $f(\tilde x_2)$ is convex over $(-\infty,\infty)$,
so its epigraph $\mathbf{epi}(f)=\mc{R}_{qmod}^{(-\infty,\infty)}(\tilde x_1-f)$
which can be represented as
\[
\left\{(\tilde x_1, \tilde x_2):\,
\exists \, y_2, z_0, z_1,
\tilde x_1 \geq y_2 + z_0,
\bbm
z_0 & z_1  & 1- z_0 \\
z_1 & 1-z_0 & \tilde x_2 -z_1 \\
1-z_0 & \tilde x_2 - z_1 & y_2 - 1 + z_0
\ebm \succeq 0
\right\}.
\]
A direct SDP representation of $S$ can be obtained by
applying the inverse perspective transformation
\[
\left\{(x_1, x_2):\,
\exists \, y_2, z_0, z_1, \,
1 \geq y_2 + z_0,
\bbm
z_0 & z_1  & x_1- z_0 \\
z_1 & x_1-z_0 & x_2 -z_1 \\
x_1-z_0 & x_2 - z_1 & y_2 - x_1 + z_0
\ebm \succeq 0
\right\}.
\]
In the above, $y_2 + z_0$ can be placed by one parameter.
The plot of the above
coincides with the shaded area in Figure~\ref{fig:bean}.
\end{example}

\subsection{General case}

For the general case, we have the following result by applying Theorem~\ref{thm:RatRep}.

\begin{theorem}
Assume $\tilde{\mc{D}}$ and $S_{\tilde{\mc{D}}}(\tilde f)$ are both convex and have nonempty interior,
%$\dim( \mc{Z}(f) \cap \{\tilde x_1=0\} \cap \pt S_{\tilde{\mc{D}}}(\tilde f)) < n-1$ and
and $\dim( \{\tilde x_1=0\} \cap \pt S_{\tilde{\mc{D}}}(\tilde f)) < n-1$.
\bdes
\item [(i)] If $\tilde f(\tilde x)$ and every $\tilde g_i(\tilde x)$
are q-module concave over $\tilde{\mc{D}}$ with respect to $\tilde x_1^{d-b}$, 
then $S_{\tilde{\mc{D}}}(\tilde f ) = \mc{R}_{qmod}^{\tilde{\mc{D}}}(\tilde f)$.

\item [(ii)] If $\tilde f(\tilde x)$ and every $\tilde g_i(\tilde x)$
are preordering concave over $\tilde{\mc{D}}$ with respect to $\tilde x_1^{d-b}$, 
then $S_{\tilde{\mc{D}}}(\tilde f) = \mc{R}_{pre}^{\tilde{\mc{D}}}(\tilde f)$.
\edes
\end{theorem}

After one perspective transformation $\mathfrak{p}$,
the singular point in $S_{\mc{D}}(f)$
is mapped to one point at infinity of $S_{\tilde{\mc{D}}}(\tilde f)$,
i.e., $S_{\tilde{\mc{D}}}(\tilde f)$ itself does not have a point
which is the image $\mathfrak{p}(0)$.
And the mapping $\mathfrak{p}$ is smooth when $x_1 >0$.
At any point $x\in S_{\mc{D}}(f)$ with $x_1 >0$,
the mapping $\mathfrak{p}$ will preserve the singularity or nonsingularity at $x$.
In this sense, the perspective transformation $\mathfrak{p}$
will remove one or more singular points.
Of course, the new convex set $S_{\tilde{\mc{D}}}(\tilde f)$
might have singularity somewhere else.
In this case, we can apply some coordinate transformation
to shift one singular point to the origin
and then apply the perspective transformation again.
So a sequence of perspective transformations might be applied.
If there are finitely many singular points on the boundary,
a finite number of perspective transformations 
can be applied to remove all the singularities.
However, this approach might not work if there are infinitely many singular points,
i.e., the singular locus is positively dimensional.
For instance, the convex set
\[
\{x \in \re^3: \,  (1-(x_1-1)^2-x_2^2)^3 - x_3^4  \geq 0 \},
\]
has a singular locus of dimension one.
In this case, a finite number of perspective transformations
is usually not able to remove all the singularies.

\section{Some discussions}

We conclude this paper with some discussions and open questions.

\medskip\noindent
{\it More general convex sets}\,
It is very natural to consider general convex sets of the form
\[
S_\mc{D}(f_1,\cdots,f_k) = \{x \in \mc{D}:\, f_1(x) \geq 0,\cdots, f_k(x)\geq 0\},
\]
where $f_1(x), \ldots, f_k(x)$ are given polynomials or rational functions
concave over the convex domain
$\mc{D} = \{x \in \re^n:\, g_1(x) \geq 0,\cdots, g_m(x)\geq 0\}$
defined by polynomials $g_1(x),\cdots, g_m(x)$. Note that 
\[
S_\mc{D}(f_1,\cdots,f_k) = S_\mc{D}(f_1) \cap \cdots \cap S_\mc{D}(f_k).
\]
So it suffices to consider each individual $S_\mc{D}(f_i)$ separately.

One interesting but unaddressed case is that 
the defining polynomials $f_i$
are concave in some neighborhood of $S_\mc{D}(f_1,\cdots,f_k)$
but neither q-module nor prepordering concave over the domain $\mc{D}$.
In this situation, is it always possible to find another domain 
$\mc{D}^\prm \supset S_\mc{D}(f_1,\cdots,f_k)$ such that
the $f_i$ is q-module or prepordering concave over $\mc{D}^\prm$
with respect to some other $(p^\prm,q^\prm)$?
Or is it always possible to find a different set of defining polynomials for 
$\mc{D} = \{x \in \re^n:\,\hat g_1(x) \geq 0,\cdots, \hat g_{m^\prm}(x)\geq 0\}$
such that $f_i$ is q-module or prepordering concave over $\mc{D}$
using new defining polynomials with respect to some different $(\hat p, \hat q)$?
%A simple and useful case in practice 
%is that $\mc{D}$ is a polyhedra.
%When $S_\mc{D}(f_1,\cdots,f_k)$ is compact convex, 
%one possible approach is to choose 
%$\mc{D}^\prm$ to be a small polyhedra containing $S_\mc{D}(f_1,\cdots,f_k)$
%so that $f_i$ is q-module or prepordering concave over $\mc{D}^\prm$
%(with respect to some polynomial pair $(p^\prm,q^\prm)$).
This is an interesting future research topic.

\medskip\noindent
{\it The separability in Positivestellensatz}\,
The rational function $f(x)$ is concave over 
$\mc{D}$ if and only if
\[
f(u) + \nabla f(u)^T(x-u) -f(x) \geq 0, \, \forall\, x, u \in \mc{D}.
\]
By Positivestellensatz of Stengle \cite{Sten}, the above is true if and only if
\begin{align*}
\eta(x,u) \cdot f_{den}(x) f_{den}^2(u) \cdot
(f(u)+\nabla f(u)^T(x-u) - f(x)) =  \hspace{3cm}  \\
\sum_{\nu\in \{0,1\}^m } g_1^{\nu_1}(x)\cdots g_m^{\nu_m}(x)
\left(  \sum_{\mu\in \{0,1\}^m } g_1^{\mu_1}(u) \cdots g_m^{\mu_m}(u)  \sig_{\nu,\mu}(x,u)  \right)
\end{align*}
for some sos polynomials $\eta(x,u), \sig_{\nu,\mu}(x,u)$.
Here $f_{den}$ is the denominator of $f(x)$ which is nonnegative over $\mc{D}$.
When $\eta(x,u) = \eta_1(x) \eta_2(u)$ is separable,
we can choose $p(x)=\eta_1(x) f_{den}(x)$ and $q = \eta_2(u) f_{den}^2(u)$,
and then get an SDP representation for $S_\mc{D}(f)$ 
by following the approach in Section~\ref{sec:RatCvx}.
However, in general case, is it always possible to find 
a factor $\eta(x,u)$ that is separable?
Or what conditions make the factor $\eta(x,u)$ to be separable?
This is an interesting future research topic.

\medskip\noindent
{\it  Resolution of singularities}\,
In algebraic geometry \cite{Har}, 
a well known result is that any singular algebraic variety
(over a ground field with characteristic zero)
is birational to a nonsingular algebraic variety.
But the convexity might not be preserved
by this birational transformation.
Given a convex semialgebraic set in $\re^n$
with singular boundary, is it is birational to 
a convex semialgebraic set with nonsingular boundary?
Or is every convex semialgebraic set in $\re^n$
equal to the projection of some higher dimensional 
convex semialgebraic set with nonsingular boundary?
To the best knowledge of the author, 
all such questions are open.
An interesting future work is to discuss how to remove 
the singular locus of convex semilagebraic sets while preserving the convexity.

\bigskip \bigskip
\noindent
{\bf Acknowledgement} \,
The author would like to thank Bill Helton for fruitful discussions.

\end{document}